\documentclass{tran-l}
\usepackage{amsfonts,amsmath,amssymb,latexsym,epsfig,graphicx}
\newtheorem{theorem}{Theorem}[section]
\newtheorem{lemma}[theorem]{Lemma}
\newtheorem{proposition}[theorem]{Proposition}
\newtheorem{corollary}[theorem]{Corollary}

\theoremstyle{definition}
\newtheorem{definition}[theorem]{Definition}
\newtheorem{example}[theorem]{Example}

\theoremstyle{remark}
\newtheorem{remark}[theorem]{Remark}

\numberwithin{equation}{section}
\title[Schubert Varieties in the Symplectic Grassmannian]{Hilbert Functions of Points on Schubert 
Varieties in the Symplectic Grassmannian}

\author{Sudhir R.~Ghorpade} 
\address{Department of Mathematics, % \newline  \indent
Indian Institute of Technology Bombay, Powai, \newline \indent
 Mumbai 400\,076, India.}
\email{srg@math.iitb.ac.in}

\author{K.~N.~Raghavan}
\address{Institute of Mathematical Sciences,
C.I.T. Campus, Taramani,   \newline  \indent
Chennai 600\,113, India.}    
\email{knr@imsc.ernet.in}


\date{September 20, 2004}
\newcommand{\ignore}[1]{}
\newcommand{\version}{\public}
\newcommand{\public}[2]{#1}
\newcommand{\mis}{{\mathfrak M}_d(V)}
\newcommand{\grass}{G_{d,2d}}
\newcommand{\gdn}{\mis}
\newcommand{\idn}{I(d)}

\newcommand{\diag}{\mathfrak{\Delta}}
\newcommand{\diagv}{\diag^v}
\newcommand{\invol}{\#}
\newcommand{\up}{\text{up}}
\newcommand{\down}{\text{down}}
\newcommand{\hash}{\invol}
\newcommand{\bhash}{\block^\invol}

\newcommand{\st}{{\,|\,}}		%general
\renewcommand{\frak}{\mathfrak}		%general
\newcommand{\surjection}{\twoheadrightarrow}%general
\newcommand{\svw}{S^v_w}		%theorem.tex
\newcommand{\svwm}{\svw(m)}		%theorem.tex
\newcommand{\pos}{{\frak N}}		%theorem.tex
\newcommand{\smvw}{\begin{mbox}
		{\rm SM}\end{mbox}^v_w}	%smt.tex
\newcommand{\smvwm}{\smvw(m)}		%smt.tex
\newcommand{\yvw}{Y_w^v}		%smt.tex
\newcommand{\ap}{{\mathfrak w}}		%proof.tex-2
\newcommand{\aptop}[1]{{\mathfrak{top}}(#1)}	%proof.tex-2
\newcommand{\apbottom}[1]{{\mathfrak{bot}}(#1)}	%proof.tex-2
\newcommand{\mon}{\frak S}		%proof.tex
\newcommand{\monhash}{{\frak S}^\hash}		%reducetwo.tex
\newcommand{\monm}{\frak M}		%reducetwo.tex
\newcommand{\monn}{\frak N}		%reducetwo.tex
\renewcommand{\part}{\frak B}		%proof.tex

\newcommand{\tv}{T^v}			%proof.tex
\newcommand{\tvwm}{T^v_w(m)}		%proof.tex
\newcommand{\smvtopv}{\begin{mbox}
		{\rm SM}\end{mbox}^{v,v}} %proof.tex	
\newcommand{\smvtopvw}{\begin{mbox}
		{\rm SM}\end{mbox}^{v,v}_w} %proof.tex	
\newcommand{\smvtopvwm}{\begin{mbox}
		{\rm SM}\end{mbox}^{v,v}_w(m)} %proof.tex	
\newcommand{\roots}{\frak R }		%proof.tex
\newcommand{\spec}{\mathfrak E}
\renewcommand{\min}{\frak T}		%proof.tex
\newcommand{\monw}{{\mon}_w}		%phi.tex
\newcommand{\piece}{\frak P}			%phi.tex
\newcommand{\blockc}{\blc}			%pi.tex
\newcommand{\block}{\frak B}		%draft.tex
\newcommand{\blc}{\frak C}		%draft.tex
\newcommand{\init}[1]{\mbox{{\rm in}}(#1)}%groebner.tex
\newcommand{\po}{>}		%groebner.tex
\newcommand{\poone}{>_1}		%groebner.tex
\newcommand{\potwo}{>_2}		%groebner.tex
\newcommand{\pothree}{>_3}		%groebner.tex
\newcommand{\pofour}{>_4}		%groebner.tex
\newcommand{\path}{\Lambda}		%interpretation.tex
\newcommand{\bstart}{\beta({\mbox{\begin{rm}start\end{rm}}})}
\newcommand{\bstartj}{\beta_j({\mbox{\begin{rm}start\end{rm}}})}
\newcommand{\bfinish}{\beta({\mbox{\begin{rm}finish\end{rm}}})}
\newcommand{\bfinishj}{\beta_j({\mbox{\begin{rm}finish\end{rm}}})}
\newcommand{\bthm}{\begin{theorem}}
\newcommand{\ethm}{\end{theorem}}

\newcommand{\blem}{\begin{lemma}}
\newcommand{\blemma}{\begin{lemma}}
\newcommand{\elemma}{\end{lemma}}
\newcommand{\elem}{\end{lemma}}

\newcommand{\bpr}{\begin{proposition}}
\newcommand{\bprop}{\begin{proposition}}
\newcommand{\epr}{\end{proposition}}
\newcommand{\eprop}{\end{proposition}}

\newcommand{\bdefn}{\begin{definition}\begin{rm}}
\newcommand{\bdefinition}{\begin{definition}\begin{rm}}
\newcommand{\edefinition}{\end{rm}\end{definition}}
\newcommand{\edefn}{\end{rm}\end{definition}}

\newcommand{\bexample}{\begin{example}\begin{rm}}
\newcommand{\eexample}{\end{rm}\end{example}}

\newcommand{\bremark}{\begin{remark}\begin{rm}}
\newcommand{\eremark}{\end{rm}\end{remark}}

\newcommand{\bcor}{\begin{corollary}}
\newcommand{\ecor}{\end{corollary}}

\newtheorem{exercise}[theorem]{Exercise}
\newcommand{\bex}{\begin{exercise}\begin{rm}}
\newcommand{\eex}{\end{rm}\end{exercise}}

\newtheorem{notation}[theorem]{Notation}
\newcommand{\bnotation}{\begin{notation}\begin{rm}}
\newcommand{\enotation}{\end{rm}\end{notation}}

\newenvironment{proofone}[1][]{\noindent{\sc Proof #1}:\ }{\ $\Box$}
\newcommand{\bpf}{\begin{proofone}}
\newcommand{\bproof}{\begin{proof}}
\newcommand{\epf}{\end{proofone}}
\newcommand{\eproof}{\end{proof}}

\newenvironment{prooftwo}[1][]{{\sc Proof #1}:\ }{}
\newcommand{\bpftwo}{\begin{prooftwo}}
\newcommand{\bprooftwo}{\begin{prooftwo}}
\newcommand{\epftwo}{\end{prooftwo}}
\newcommand{\eprooftwo}{\end{prooftwo}}

\newenvironment{solution}{{\sc Solution}:\ }{\ $\Box$}
\newcommand{\bsol}{\begin{solution}}
\newcommand{\esol}{\end{solution}}

\begin{document}

\begin{abstract}
We give an explicit combinatorial description of the multiplicity as well as 
the Hilbert function of the tangent cone at any point on a Schubert variety 
in the symplectic Grassmannian. 
\end{abstract}

\maketitle
%\baselineskip=17pt
%\input{intro.tex} 

%\begin{abstract}
%We give an explicit combinatorial description of the multiplicity as well as 
%the Hilbert function of the tangent cone at any point on a Schubert variety 
%in the symplectic Grassmannian. 
%\end{abstract}

\section{Introduction} 
\label{sintro}

Let $G$ be a semisimple algebraic group over an algebraically closed field 
$k$ and $P$ be a parabolic subgroup of $G$. Fix a Borel subgroup 
$B$ of $G$ and a
maximal torus $T$ in $G$ such that $T\subset B\subset P$. Now, $G/P$ is 
projective variety and it has a distinguished class of subvarieties, known as
Schubert varieties (in $G/P$); these are indexed by a set $W^P$, which 
corresponds to the $T$-fixed points of $G/P$ for the action given by left
multiplication. Given $w\in W^P$, we denote by $e^w$ the corresponding 
$T$-fixed point, and by $X_w$ the corresponding Schubert variety. In fact, 
$X_w$ is the closure of the $B$-orbit $Be^w$, and can be decomposed as the
union of `smaller' $B$-orbits $Be^v$, where $v$ ranges over elements of 
$W^P$ satisfying $v\le w$; here, $\le$ is a certain partial order, called
the Bruhat-Chevalley order. 

The study of Schubert varieties, in general, and the singularities of 
Schubert varieties, in particular, has been an active and vibrant area
of research in the past three decades.  It may suffice to cite the recent 
monograph \cite{BL} by Billey and Lakshmibai, which surveys many known 
results and can also be a useful reference for the background material.
Among the basic questions, insofar as the singularities of Schubert varieties 
are concerned, are the following:  (1) Which points are singular? (2) what 
is the multiplicity at a (singular) point? and (3) what is the Hilbert 
function (of the tangent cone) at a (singular) point? It may be noted that 
these questions are in an
ascending order of generality since the singular points are those 
of multiplicity $>1$, and the Hilbert function determines the Hilbert
polynomial whose (normalized) leading coefficient gives the multiplicity.
Also note that, thanks to the $B$-orbit decomposition, it suffices to only
look at the points $e^v$ in $X_w$ where $v\le w$. 

The singular loci of Schubert varieties are fairly well understood, 
thanks to the works of several mathematicians 
(see \cite{BL} for details). Recursive formulas for the multiplicity and for 
the Hilbert function in the case of minuscule $G/P$ and also in the case of 
symplectic Grassmannian were obtained by Lakshmibai and Weyman \cite{lw} in 
1990. The singularities of Schubert varieties in the symplectic Grassmannian have also
been studied by Brion and Polo \cite{bp} who determine the multiplicity at `generic'
singular points.  More explicit results seem to be available so far in the special case
when $G={\rm SL}_n = {\rm SL}(V)$, where $V$ is an $n$-dimensional vector 
space over $k$ and $P=P_d$ the maximal parabolic subgroup 
given by those $g\in {\rm SL}(V)$ which stabilize a (fixed) $d$-dimensional 
subspace of $V$. Assume for a moment that we are in this case. Thus $G/P$ 
is the Grassmannian $G_{d}(V)$ and $W^P$ may be viewed as the set $I(d,n)$ 
of subsets of $\{1, \dots , n\}$ of cardinality $d$. Let $\epsilon$ denote 
the element $\{1, \dots, d\}$ of $I(d,n)$. 
An explicit closed-form formula for the multiplicity $m_{\epsilon}(w)$ 
of $X_w$ at $e^{\epsilon}$ was given by Lakshmibai and Weyman  \cite{lw}
in their 1990 paper. 
Recently, an explicit description for the Hilbert function of 
$X_w$ at $e^{\epsilon}$ 
was given by Kreiman and Lakshmibai \cite{kl}. Further, in the case of 
arbitrary $v$, Kreiman and Lakshmibai formulated 
two conjectures which give an explicit combinatorial description for
(i)  the Hilbert function of  $X_w$ at $e^v$ and
(ii)  the multiplicity $m_{\epsilon}(w)$ of  $X_w$ at $e^v$.
Subsequently, both the conjectures were
proved by Kodiyalam and Raghavan \cite{kr}, and independently by Kreiman \cite{kreiman} 
(see also \cite{klr}).  Both \cite{kr} and \cite{kreiman} also obtain a reformulation
of the main result in terms of Gr\"obner bases. 

In this paper, we consider the problem of determining the multiplicity
as well as the Hilbert function of Schubert varieties in the symplectic
Grassmannian $G/P$ where $G$ is the symplectic group ${\rm Sp}_{2d}$
and $P$ a maximal parabolic subgroup $P_d$ of $G$. 
Alternatively, $G/P$ is the Grassmannian of maximal isotropic subspaces. 
[see Section \ref{stheorem} for details]. 
For Schubert varieties in symplectic Grassmannians, we give
an explicit combinatorial description for the multiplicity, and more 
generally, the Hilbert function at any $e^v$ in $X_w$. In effect, we 
formulate and prove the two conjectures of Kreiman and Lakshmibai in this
case.  Moreover, as in \cite{kr}, we obtain a reformulation of the main 
result in terms of Gr\"obner bases. 
A precise statement of our main result is given in Section 
\ref{stheorem}. 

A key ingredient in our proof is the {\em Standard Monomial Theory} 
for symplectic Grassmannians and its Schubert subvarieties, as developed 
by  De~Concini~\cite{dc} and by Lakshmibai, Musili and Seshadri~\cite{lms}. 
This allows us to translate the problem from geometry to combinatorics.   This
translation---although we derive it afresh in Section \ref{ssmt} below for the 
sake of clarity, completeness, and
readability---is already there in Lakshmibai-Weyman~\cite{lw}.
Our main job then is the solution of the resulting combinatorial
problem.  Roughly speaking,  it amounts to showing that the
combinatorial constructions of \cite{kr} behave well with respect 
to a certain involution.   
This is carried out in Sections \ref{sfurtherred} and \ref{sproofcompletion}.
Moreover, as in \cite{kl} and \cite{kr}, we obtain an interpretation of 
the multiplicity as the number of certain nonintersecting lattice paths, and this
is described in Section \ref{sinterpretations}. 

The interpretation of the multiplicity in terms of nonintersecting lattice paths 
may be viewed as an analogue of the results of Krattenthaler \cite{Kra1,Kra2}.
To review the latter, we note that in the case of the classical Grassmannian,  
Rosenthal and Zelevinsky \cite{rz} obtained in 1998 a closed-form formula for 
$m_v(w)$ for arbitrary $v$ and $w$, using the recursive formula of Lakshmibai and Weyman  \cite{lw}. The lattice path interpretation was used by Krattenthaler \cite{Kra1} 
to explain the relationship between the Rosenthal--Zelevinsky formula and the 
Lakshmibai-Weyman formula for $m_{\epsilon}(w)$. Further, 
Krattenthaler \cite{Kra2} also used the lattice path interpretation to 
prove the Multiplicity Conjecture of Kreiman and Lakshmibai. 

We now attempt to outline a connection of Schubert varieties, in general, and the results
of this paper, in particular, to a class of affine varieties, broadly known as determinantal
varieties. For the last two decades, the study of determinantal
varieties has proceeded almost in parallel to the study of 
Schubert varieties, but often as an independent pursuit. 
In particular, explicit formulas for the multiplicity and the Hilbert 
function for various classes of determinantal varieties have been obtained. 
The relation between determinantal varieties
and Schubert varieties is best explained in the case of Grassmannian
$G_d(V) = SL_n/P$. We have the well-known Pl\"ucker embedding of $G_d(V)$ 
in ${\mathbb P}(\bigwedge^dV)$ and for each $v\in W^P$, there is a basic 
affine open set ${\mathbb A}_v$ containing $e^v$. In case, $v=\epsilon$,
the intersections $X_w \cap {\mathbb A}_v$ are affine varieties, which are
precisely the varieties defined by an ideal `cogenerated' by a minor of a
generic matrix; 
in other words, these are exactly the determinantal varieties studied 
by Abhyankar \cite{abh}, and later by Herzog-Trung \cite{ht}, 
Conca-Herzog \cite{CoH}, and others (see, e.g., \cite{Gh1,Gh2}).
Moreover, this affine variety is a cone and hence it coincides with the 
tangent cone to $X_w$ at $e^{\epsilon}$. 
Thus, in retrospect, the multiplicity formula of Lakshmibai-Weyman and the
Hilbert function result of Kreiman-Lakshmibai could have been deduced from 
the work of Abhyankar \cite{abh} and others. On the other hand, for arbitrary
$v$, the varieties $X_w \cap {\mathbb A}_v$ are not so well understood and
as far as we know, there is no analogue in the literature on determinantal 
varieties of the results in \cite{kr}. As remarked in \cite{kr}, 
the recent work of Knutson and Miller \cite{km} considers a class
of determinantal varieties more general than those cogenerated by a minor, but
it is not clear if the varieties $X_w \cap {\mathbb A}_v$ belong to this class
when $v\ne \epsilon$. 

In a similar vein, considering affine patches
of Schubert varieties in the symplectic Grassmannian, when $v=\epsilon$, 
leads to varieties given by ideals cogenerated by a minor of a generic 
symmetric matrix. Thus in this special case, the results obtained in this
paper could be compared with those of Conca \cite{conca} on symmetric 
determinantal varieties.  Likewise, one could
take up the case of orthogonal Grassmannians (although we do not do this here)
and the results thus obtained ought to be compared with those of Ghorpade and 
Krattenthaler \cite{GK} on pfaffian varieties. In either of these cases, for
an arbitrary $v$, the varieties $X_w \cap {\mathbb A}_v$ do not seem to 
correspond to any of the known classes of varieties defined by the minors of
a generic symmetric matrix or the pfaffians of a generic skew-symmetric matrix. 

There is a yet another related, but independent, body of work on 
degeneracy loci. Formulas for the fundamental classes of degeneracy
loci of maps of vector bundles give rise to multiplicity formulas for
determinantal and pfaffian varieties. For a detailed explanation, we
refer to the appendix in \cite{GK} and the books of Fulton-Pragacz \cite{FP}
and Manivel \cite{Ma}. While some of the results on multiplicity could be 
deduced from the corresponding results on degeneracy loci, it does
not seem likely that the latter impinge on the determination of 
Hilbert functions.

\section{The Theorem}\label{stheorem}
The main results of this paper are stated in Theorem~\ref{tmain} and
Corollary~\ref{cmain} below. But first, we fix some notation and 
terminology to be used in the rest of this paper and briefly review some 
preliminary notions and results.
% provide respectively combinatorial
%expressions for the Hilbert function and the multiplicity at a point
%on a Schubert subvariety of the variety $\gdn$ of maximal isotropic
%subspaces.   In order to state them, we start by recalling some well known facts about Schubert
%varieties in $\gdn$.    This also serves to fix notation.

%We keep the notations of \S1.  In other words, \ldots
Given any nonnegative integer $n$, we denote by $[n]$ the set $\{1,\dots ,n\}$. 
The cardinality of a finite set $v$ will be denoted by $|v|$. Given positive
integers $r$ and $n$ with $r\le n$, we denote by $I(r, n)$ the set of all subsets 
of $[n]$ of cardinality $r$. An element $v$ of $I(r,n)$ may be written as 
$v=(v_1,\ldots,v_r)$ where $1\leq v_1<\ldots<v_r\leq n$ and $v= \{v_1, \dots , v_r\}$. 
Given any $v= (v_1,\ldots ,v_r)$ and $w=(w_1, \ldots ,w_r)$ in $I(r,n)$, we say 
$v\leq w$ if $v_1\leq w_1$, $\ldots$, $v_r\leq w_r$. Clearly, $\le$ defines a partial 
order on $I(r,n)$.

A positive integer $d$ will be kept fixed throughout this paper. 
For $j \in [2d]$,  set $j^*:=2d+1-j$.
Let $\idn$ denote the set of subsets $v$ of $[2d]$ with the
property that exactly one of $j$,~$j^*$ belongs to $v$ for
every $j\in [d]$.  Note that $I(d)$ is a subset of $I(d,2d)$. 
In particular, we have the partial order $\le$ on $I(d)$ induced from $I(d,2d)$. 
We denote by $\epsilon$ the element $[d]=(1, \dots, d)$ of $I(d)$. 

Fix a vector space $V$ of dimension $2d$ over an algebraically closed
field of arbitrary characteristic.  Fix a nondegenerate skew-symmetric
bilinear form $\langle\ ,\ \rangle$ on $V$.   Fix a basis $e_1, \ldots, e_{2d}$
of $V$ such that 
\[ 
\langle e_i, e_j \rangle = \left\{ \begin{array}{rl}
	1 & \mbox{ if $i=j^*$ and $i<j$}\\
	-1 & \mbox{ if $i=j^*$ and $i>j$}\\
	0 & \mbox{ otherwise}\\ \end{array}\right.
 \]
A linear subspace of $V$ is said to be {\em isotropic} if the form  
$\langle\ ,\ \rangle$ vanishes identically on it. It is well-known that 
an isotropic subspace of $V$ has dimension at most $d$ and every isotropic 
subspace is contained in one of dimension $d$.   Denote by $G_d(V)$ the 
Grassmannian of $d$-dimensional subspaces of $V$ and by $\mis$ 
the set of all $d$-dimensional isotropic subspaces of $V$. Then $\mis$ is a 
closed subvariety of $G_d(V)$, and is called the {\em variety of maximal 
isotropic subspaces} or the {\em symplectic Grassmannian}.
  
The group ${\rm Sp}(V)$ of linear automorphisms of $V$ preserving
$\langle\ ,\ \rangle$ acts transitively on $\mis$---this follows from
Witt's theorem that an isometry between subspaces can be lifted to one
of the whole vector space.  So $\mis$ is identified as the quotient of 
${\rm Sp}(V)$ by the stabilizer of any point (for example, the span of 
$e_1, \ldots, e_d$).   The elements of ${\rm Sp}(V)$ that are diagonal 
with respect to the basis $e_1$, \ldots, $e_{2d}$ form a maximal torus 
$T$ of ${\rm Sp}(V)$.   Similarly the elements of ${\rm Sp}(V)$ that are 
upper triangular with respect to %the basis 
$e_1$, \ldots, $e_{2d}$ form a Borel subgroup $B$ of ${\rm Sp}(V)$---a 
linear transformation is {\em upper triangular} with respect to 
$e_1,\ldots,e_{2d}$ if for each $j\in [2d]$,  the image of $e_j$  under 
the transformation is a linear combination of %$e_1, \dots , e_j$.
$e_i$ with $i\leq j$.

The $T$-fixed points of $\mis$ are parametrized by $I(d)$: for
$v=(v_1,\ldots,v_d)$ in $I(d)$,  the corresponding $T$-fixed point,
denoted by $e^v$,  is the span of $e_{v_1},\ldots,e_{v_d}$.
These points lie in different $B$-orbits and the union of 
their $B$-orbits %the $B$-orbits of these points 
is all of $\gdn$.
A {\em Schubert variety} in $\gdn$ is by definition the closure of
such a $B$-orbit (with the reduced scheme structure).     
Schubert varieties are thus indexed by the $T$-fixed points
and so in turn by $\idn$.  Given $w$ in $\idn$,  we denote by $X_w$ the closure
of the $B$-orbit of the $T$-fixed point $e^w$. We have the $B$-orbit decomposition:
$$
X_w = \coprod_{v\le w} Be^v.
$$
We are interested in the local rings of various points on a Schubert
variety $X_w$.    In view of the above $B$-orbit decomposition, it is enough to 
focus attention on the $T$-fixed points contained in $X_w$, that is, the points
$e^v$ for $v\le w$.

For the rest of this section,  fix elements $v,w$ of $\idn$ with
$v\leq w$. Define 
$$
\roots^v:=\left\{ (r,c)\in  [2d]\setminus v \times v \, : \, r\le c^*\right\} \quad
{\rm and} \quad 
\pos^v:=\left\{ (r,c)\in  [2d]\setminus v \times v \, : \, r > c\right\}.
$$
We will be considering ``multisets'' on $\roots^v$ and $\pos^v$.
By a {\em multiset} on a finite set $S$ we mean a  collection of elements of $S$ 
in which repetitions are allowed and kept account of.   Multi-sets on $S$ can be thought 
of as monomials in the variables corresponding to the elements of $S$. The cardinality 
of a multiset is the number of elements in it, counting repetitions, or equivalently,
the degree of the corresponding monomial.  The union of multisets is the product
of the corresponding monomials.   The intersection of a multiset with a subset
is again easily described in terms of monomials:  set equal to $1$
those variables that do
not belong to the subset.

Given any $\beta_1=(r_1,c_1)$, $\beta_2=(r_2,c_2)$ in $\pos^v$,
we say that $\beta_1>\beta_2$ if $r_1>r_2$ and $c_1<c_2$.   A sequence
$\beta_1>\ldots>\beta_t$ of elements of $\pos^v$ is called a
{\em $v$-chain}.
Given a $v$-chain $\beta_1=(r_1,c_1)>\ldots>\beta_t=(r_t,c_t)$, we define
$$
s_{\beta_1}\cdots s_{\beta_t}v := \left( \{v_1,\ldots,v_d\}\setminus\{c_1,\ldots,c_t\}\right)\cup\{r_1,\ldots,r_t\},
$$
and note that this is an element of $\idn$.  
In case the $v$-chain is empty,  this element is just $v$.
We say that $w$ {\em dominates} the $v$-chain
$\beta_1>\ldots>\beta_t$  if $w\geq 
s_{\beta_1}\cdots s_{\beta_t}v$.

Let $\mon$ be a monomial on $\roots^v$.  By a {\em $v$-chain in $\mon$}
we mean a sequence $\beta_1>\ldots>\beta_t$ of elements of
$\mon \cap\pos^v$. 
We say that {\em $w$ dominates $\mon$} if $w$ dominates every $v$-chain in $\mon$.

Let $S^v_w$ denote the set of $w$-dominated monomials on $\roots^v$,
and $S^v_w(m)$ the set of such monomials of degree $m$.

We can now state our theorem:
\bthm \label{tmain} 
Let $v,w$ be elements of $\idn$ with $v\leq w$.   Let $X_w$ be the
Schubert variety corresponding to $w$,  $e^v$ the $T$-fixed point
corresponding to $v$, and $R$ the coordinate ring of the tangent cone
to $X_w$ at the point $e^v$ (that is, the associated graded ring
${\begin{rm}gr\end{rm}}_{{M}}({\mathcal O}_{X_w, e^v})$ of the local ring
${\mathcal O}_{X_w,e^v}$ of $X_w$ at the point $e^v$ with respect to its 
maximal ideal $M$).  Then the dimension as a vector space of the 
$m^{\begin{rm}th\end{rm}}$ graded piece $R(m)$ of $R$ equals the 
cardinality of $\svwm$, where $\svwm$ is as defined above.
\ethm

The proof of the theorem occupies
sections~\ref{ssmt},~\ref{sfurtherred},~and~\ref{sproofcompletion}.   
For now let us note the following easy consequence.
\bcor\label{cmain}\label{ctmain}
With notation as in Theorem~\ref{tmain} above,
the multiplicity of $R$ equals the number of square-free
$w$-dominated monomials on $\roots^v$ of maximum cardinality.
\ecor
\bproof
The proof of the corresponding corollary in \cite{kr}
works verbatim here.
\eproof

\section{Reduction to combinatorics}
\label{ssmt}

Let $\mis\subseteq G_{d}(V)\hookrightarrow \mathbb{P}(\wedge^d V)$ be 
the Pl\"ucker embedding.  
The homogeneous coordinate rings of $\mis$ and its Schubert subvarieties
in this embedding have been described by De~Concini~\cite{dc} and 
Lakshmibai, Musili, and Seshadri~\cite{lms}.
We will use their results to reduce the proof of Theorem~\ref{tmain} to
combinatorics.   Our primary reference will be \cite{dc}---its language
and approach suit our purpose well.

For $\theta$ in $I(d,2d)$,  let $p_\theta$ denote the corresponding Pl\"ucker
coordinate. Consider the affine patch $\mathbb{A}$ of $\mathbb{P}(\wedge^d V)$
given by $p_{\epsilon}=1$.     The intersection $\mathbb{A}\cap\grass$ of
this patch with the Grassmannian is an affine space.     Indeed the $d$-plane
corresponding to an arbitrary point $z$ of $\mathbb{A}\cap\grass$ has a basis consisting
of column vectors of a matrix of the form
\[ B= \left(\begin{array}{c}I\\A\\ \end{array}\right)\]
where $I$ is the identity matrix of size $d\times d$  and
$A$ is an arbitrary matrix of size $d\times d$.
The association $z\mapsto A$ is bijective.   
The restriction
of a Pl\"ucker coordinate $p_\theta$ to $\mathbb{A}\cap\grass$ is given by
the determinant of a submatrix of size $d\times d$ of $B$,  the entries of 
$\theta$ determining the rows to be chosen from $B$ to form the submatrix.

As can be readily verified, a point $z$ of $\mathbb{A}\cap\grass$ belongs to 
$\mis$ if and only if the corresponding matrix $A=(a_{ij})$ is {\em symmetric 
with respect to the anti-diagonal\/}:   $a_{ij}=a_{j^*i^*}$,  where the 
columns and rows of $A$ are numbered $1,\ldots,d$ and $d+1,\ldots,2d$ 
respectively, and $r^*:=2d-r+1$ for $r\in [2d]$. 
For example, if $d=4$, then a matrix that is symmetric with respect to the
anti-diagonal looks like this:
\[
\left(\begin{array}{cccc}
d & c & b & a\\
g & f & e & b\\
i & h & f & c\\
j & i & g & d\\ \end{array}\right)	\]
      
\bnotation 
For $u=(u_1,\ldots,u_d)$ in $I(d,2d)$,  set $u^*:=(u_d^*,\ldots,u_1^*)$.
The association $u\mapsto u^*$ is an order reversing involution of $I(d,2d)$.
There is another order reversing involution on $I(d,2d)$, namely 
$u\mapsto [2d]\setminus u$.   These two involutions commute with
each other.   Composing them, we obtain an order preserving involution
on $I(d,2d)$:  $u\mapsto u^\hash:=[2d]\setminus u^*$. Note that $u\in I(d)$ if and
only if $u =  u^\hash$.
\enotation

\blemma
\label{lpluecker} 
Given any $\theta \in I(d)$, we have
$p_\theta=p_{\theta^\hash}$ on $\mis$.
\elemma

\bproof
Since $\mis$ is irreducible and
its intersection $\mathbb{A}\cap\mis$ with the affine patch is non-empty,  it is
enough to check that $p_\theta=p_{\theta^\hash}$ on $\mathbb{A}\cap\mis$,  and this
follows from the symmetry property just mentioned of the matrix $A$.
\eproof

The relations $p_\theta=p_{\theta^\hash}$ do not span the space of all linear
relations among the $p_\theta$---see Example~\ref{elw} below.
In order to describe a nice parametrizing set for a basis for the space of
linear forms in the 
homogeneous coordinate ring of $\mis$---in fact for describing bases
for spaces of forms of any given degree---we make the following definition.     

\begin{definition}\hspace{-1.5mm}\footnote{
Admissible pairs as defined here are
a special case of the {\em admissible minors\/} of De~Concini~\cite{dc}:
we are only considering the case $k=r$ in his notation.   Our
standard tableaux are his {\em standard symplectic tableaux\/}
but here again we are only considering the special case $k=r$.

The original definition of admissible pairs by
Lakshmibai-Musili-Seshadri~\cite[Part\,A,~\S3]{lms} is in the more general
context of a quotient by a maximal parabolic subgroup of classical type of a
semisimple algebraic group.
The realization of the importance of
admissible pairs was a key step in their development of standard monomial
theory.
The definition given here is equivalent,  in the special
case being considered,  to theirs.

In Littelmann's language of {\em paths\/}~\cite{littinv,littann},
an admissible pair is just an {\em L-S path\/} of shape a fundamental weight of classical 
type.
}
\label{dap}\begin{rm}
The {\em $\epsilon$-degree\/} of an element $x$ of $I(d)$ is the
cardinality of $x\setminus [d]$ or equivalently that of 
$[d]\setminus x$. More generally, given any $v \in I(d)$, the {\em $v$-degree\/} of 
an element $x$ of $I(d)$ is the cardinality of $x\setminus v$ or equivalently that 
of $v\setminus x$.
An ordered pair $\ap=(x,y)$ of elements of $I(d)$ is called an {\em admissible pair\/} if
$x\geq y$ and the $\epsilon$-degrees of $x$ and $y$ are equal.  We refer to
$x$ and $y$ as the {\em top\/} and the {\em bottom\/} of $\ap$ and write 
$\aptop{\ap}$ for $x$ and $\apbottom{\ap}$ for $y$. 

Given any admissible pairs
$\ap=(x,y)$ and $\ap'=(x',y')$, we say $\ap \ge \ap'$  if $y\ge x'$, that is, if
$x\ge y\ge x'\ge y'$.
An ordered sequence $(\ap_1 ,\ldots ,\ap_t)$
of admissible pairs is called a {\em standard tableau\/} if 
$\ap_i\ge \ap_{i+1}$ for $1\le i <t$.
We often write $\ap_1\geq\ldots\geq\ap_t$ to denote
the standard tableau $(\ap_1 ,\ldots ,\ap_t)$.
Given any $w\in I(d)$, we say that 
a standard tableau $\ap_1\geq\ldots\geq\ap_t$  is {\em $w$-dominated\/} 
if $w\geq \aptop{\ap_1}$.  
\edefinition

\bprop\label{pap}
There is a bijective map $(x,y)\mapsto(\theta,\tau)$ from the set of 
admissible pairs $(x,y)$ 
onto the set of ordered pairs $(\theta,\tau)$ 
of elements of $I(d,2d)$ satisfying
$$
\tau=\theta^\hash, \quad |\theta\cap [d]| = |\theta^\hash\cap [d]|,  \quad 
{\rm and} \quad \theta\cap [d] \ge \theta^\hash\cap [d].
$$
\eprop
\bproof  
Let $[d]^{\rm c}: = [2d]\setminus [d] =\{d+1, \dots , 2d\}$.
Given an admissible pair $(x,y)$, set
$$
\theta = \left( x\cap [d] \right) \cup \left( y\cap [d]^{\rm c} \right)
\quad {\rm and} \quad 
\tau = \left( y\cap [d] \right) \cup \left( x \cap [d]^{\rm c} \right)
$$
That $(\theta,\tau)$ satisfies the three conditions is readily verified.
For the map in the other direction,  set
$$
x = \left( \theta \cap [d] \right) \cup \left(\tau \cap [d]^{\rm c} \right)
\quad {\rm and} \quad 
y = \left( \tau \cap [d] \right) \cup \left( \theta \cap [d]^{\rm c} \right)
$$
It is easy to verify that $(x,y)$ is an admissible pair.\eproof

\bexample  \label{elw}
The relations $p_\theta=p_{\theta^\hash}$ do not span the linear space of
relations although \cite[page~198, item~(b)]{lw} seems to claim just that.
For example, let 
$$
d=4,  \quad 
\theta_1=(1,2,7,8),  \quad 
\theta_2=(1,4,5,8)  \quad \mbox{ and } \quad 
\theta_3=(1,3,6,8).
$$
Then 
$$
\theta_1^\hash=(3,4,5,6), \quad 
\theta_2^\hash=(2,3,6,7)  \quad \mbox{ and } \quad 
\theta_3^\hash=(2,4,5,7).
$$
{From} the form displayed above of a matrix that
is symmetric with respect to the anti-diagonal,   
we get
$$
p_{\theta_1}=p_{\theta_1^\hash}=df-cg, \quad 
p_{\theta_2}=p_{\theta_2^\hash}=cg-bi  \quad \mbox{ and } \quad 
p_{\theta_3}=p_{\theta_3^\hash}=bi-df, 
$$
and consequently,  
$p_{\theta_1} + p_{\theta_2} + p_{\theta_3} =0$.

Again, contrary to what is claimed in \cite[item~(10), page~199]{lw}, 
the association of the previous proposition is not a bijection
if the second and third conditions are dropped.
In the case $d=4$ for example there are $42$ admissible
pairs $(x,y)$ and $43$ pairs $(\theta,\theta^\hash)$.   If we try the 
procedure in the proof above for recovering $(x,y)$ on an arbitrary
$(\theta,\theta^\hash)$,
the two resulting elements may not be comparable.    Taking for example
$\theta=(2,3,6,7)$,  $\theta^\hash=(1,4,5,8)$,   we recover 
$(2,3,5,8)$ and $(1,4,6,7)$.\eexample

\bdefinition
Given an admissible pair $\ap=(x,y)$,  we define the associated Pl\"ucker coordinate
$p_\ap$ to be $p_\theta$ where $(x,y)\mapsto(\theta,\theta^\hash)$ by
the association of Proposition~\ref{pap}.   This is well-defined since
$p_\theta=p_{\theta^\hash}$ by Lemma~\ref{lpluecker}.   If $x=y$ we sometimes
write $p_x$ for $p_{(x,x)}$.  To formal products of
admissible pairs we associate the product of the associated Pl\"ucker coordinates.
In particular,  to a standard tableau $\ap_1\geq\ldots\geq\ap_t$ we associate
$p_{\ap_1}\cdots p_{\ap_t}$.    Such monomials associated to standard
tableaux are called {\em standard monomials\/}. Given any $w$ in $I(d)$, we say that 
a standard monomial is {\em $w$-dominated} if the corresponding standard
tableau is $w$-dominated.
\edefinition

We can now state the main theorem of standard monomial theory for $\mis$ and
its Schubert subvarieties.

\bthm
\label{tsmt} {\bf (De~Concini, Lakshmibai-Musili-Seshadri)}
Standard monomials
$p_{\ap_1}\cdots p_{\ap_t}$ of degree $r$ form a basis for the
space of forms of degree $r$ in
the homogeneous coordinate ring of $\mis$ in the Pl\"ucker embedding.
More generally,   given any $w \in I(d)$,  the
$w$-dominated standard monomials of degree $r$ 
form a basis for the space of forms of 
degree $r$ in the homogeneous coordinate ring of the Schubert subvariety
$X(w)$ of $\mis$.
\ethm

\bproof  
The linear independence of the $w$-dominated standard monomials in the
homogeneous coordinate ring of $X_w$ is proved 
in \cite[Lemma~3.5]{dc}\footnote{ 
The Schubert varieties in \cite{dc} are orbits under the lower triangular
Borel subgroup as opposed to our choice of upper triangular here and
that is why the domination is reversed in the statement of the Lemma
there.}.
That all standard monomials span the
homogeneous coordinate ring of $\mis$ (and so also that of $X_w$)
is the content of 
\cite[Theorem~2.4]{dc};   on the other hand,  as is easy to see,
$p_\vartheta$ vanishes on $X_w$ unless $w\geq\aptop{\vartheta}$,
and so the standard
monomials that are not $w$-dominated vanish on $X_w$.\eproof

{From} the above theorem we now deduce a basis for the coordinate ring for
an affine patch of a Schubert subvariety in $\mis$.   

\bdefinition 
Given any $v\in I(d)$, we say that a standard tableau 
$\ap_1\geq\ldots\geq\ap_t$ is {$v$-compatible \/} if for each $\ap_i$, 
either $v\geq\aptop{\ap_i}$ or $\apbottom{\ap_i}\geq v$,  and
$\ap_i\neq (v,v)$.
A standard monomial is {\em $v$-compatible \/} if the corresponding standard
tableau is $v$-compatible.   Given $v$ and $w$ in $I(d)$, we denote by $\smvw$ the set
of $w$-dominated $v$-compatible standard tableaux.
\edefinition

Fix elements
$v\leq w$ of $I(d)$ so that the point $e^v$ belongs to the 
Schubert variety $X(w)$.   Let $\mathbb{A}_v$ denote the
affine patch of $\mathbb{P}(\wedge^d V)$ given by $p_v\neq 0$ and set 
$$\yvw:=X(w)\cap\mathbb{A}_v.$$    
The point $e^v$ is the
origin of the affine space $\mathbb{A}_v$.   

The functions $f_\ap=p_\ap/p_v$, $\ap$ an admissible pair, provide a 
set of coordinate functions on $\mathbb{A}_v$.   The coordinate ring $k[\yvw]$
of $\yvw$ is a quotient of the polynomial ring $k[f_\ap]$, where $k$ is 
the underlying field.      
\bprop \label{psmt}
As $\ap_1\geq\ldots\geq\ap_t$
runs over the set ${\text SM}^v_w$ of
$w$-dominated $v$-compatible standard tableaux,
the elements
$f_{\ap_1}\cdots f_{\ap_t}$
form a basis for the coordinate
ring $k[\yvw]$ of the affine patch $Y_w^v=X_w\cap\mathbb{A}^v$ of the
Schubert variety~$X_w$.\eprop
\bproof
The proof is similar to the proof of Proposition~3.1 of \cite{kr}.
First consider any linear dependence relation among the
$f_{\ap_1}\cdots f_{\ap_t}$.     Replacing $f_\ap$ by $p_\ap$ and ``homogenizing'' by 
$p_v$ yields a linear dependence relation among the
$w$-dominated standard monomials
$p_{\vartheta_1}\cdots p_{\vartheta_s}$ restricted to $X_w$,
and so the original relation must only have been the trivial one, for
by Theorem~\ref{tsmt}
the 
$p_{\vartheta_1}\cdots p_{\vartheta_s}$  are linearly independent.

To prove that the
$f_{\ap_1}\cdots f_{\ap_t}$ span $k[Y_w^v]$ as a vector space,
we need to look at not only the corresponding statement for the Pl\"ucker
coordinates but also the proof of that statement.
What is immediate from the corresponding statement for the Pl\"ucker 
coordinates is that 
$f_{\vartheta_1}\cdots f_{\vartheta_s}$ span $k[\yvw]$ as 
$\vartheta_1\geq\ldots\geq\vartheta_s$ varies over $w$-dominated
standard monomials---the problem at hand is to show that 
the $v$-compatible ones among these are enough.

To an arbitrary monomial 
$p_{\vartheta_1}\cdots p_{\vartheta_s}$ in the Pl\"ucker coordinates,
attach the following multiset of $\{1,\ldots,2d\}$:
\[
\aptop{\vartheta_1}\cup\apbottom{\vartheta_1}\cup
\cdots\cup\aptop{\vartheta_s}\cup\apbottom{\vartheta_s}
\]
We claim that if 
$p_{\varphi_1}\cdots p_{\varphi_s}$
is a standard monomial that occurs with non-zero coefficient in the
expression of 
$p_{\vartheta_1}\cdots p_{\vartheta_s}$ as a linear combination
of standard monomials,    the multiset attached to
the two monomials are the same.
The claim follows from the
nature of the relations used in the proof in \cite{dc}
of the spanning by the standard
%%%%%%%%%% A Long Footnote %%%%%%%%%%%%%%%%%%%%%%%%%%%%%%%%%%%%%%
monomials\footnote{There are two types of relations used in the proof:  
those in 
equation~(1.1) and those in Proposition~1.8 (all numbers as in \cite{dc}).
To get the theorem about spanning by standard monomials in our situation,
we only need to use special cases of these and so let us first specialize.

In~(1.1) we take $s=s'$ and further let these equal our $d$ (our $d$
is De~Concini's~$r$);   the latter halves of all minors are $1,\ldots,d$
for us and so we write only the first half---a minor for us is therefore
just an element
$\theta$ of $I(d,2d)$;   the right side $H$ is $0$
in our case by the choice $s=d$.    Equation~(1.1) gets used in the
following way:   %the left side is a linear combination
%of products of two minors and the right side is zero;  
whenever we have a `bad' product of two minors,
it occurs as a term on the left side of an equation of type~(1.1) where
the other terms are `not so bad'; the right side being $0$,
this allows replacement of 
a `bad' product by a linear combination of `not so bad' ones.
Observe the following: for each term on the left side,
the multiset of $\{1,\ldots,2d\}$
that is the union of the indices of the
two minors, is constant for all terms.

Let us denote a minor $\theta$ by the pair $(\theta,\theta^\hash)$
rather than by just $\theta$.    For each term $\theta\tau$
on the left side of~(1.1)
let us consider the multiset $\theta\cup\theta^\hash\cup\tau\cup\tau^\hash$;
it follows from the observation in the last
line of the previous paragraph that this multiset is constant for all
terms.

We need Proposition~1.8 only in the case where $k$ is $d$
and $h_1,\ldots,h_k$ is $1,\ldots,d$
(we will omit writing the $h_1,\ldots,h_k$).  
Let $\theta$ is the element of $I(d,2d)$ that denotes
the minor that is
denoted $(\tilde{J}\cup\Gamma, \tilde{I}\cup\Gamma)$ in \cite{dc}
($\tilde{J}$, $\tilde{I}$, $\Gamma$ are subsets of $\{1,\ldots,d\}$;
$\Gamma$ does not meet $\tilde{J}\cup\tilde{I}$;   it is allowed
that $\tilde{I}\cap\tilde{J}$ is non-empty).
Write $M:=(\tilde{I}\cap\tilde{J})\cup\Gamma$,
$N:=\{1,\ldots,d\}\setminus(\tilde{I}\cup\tilde{J}\cup\Gamma)$,
$\hat{I}:=\tilde{I}\setminus \tilde{J}$,
and $\hat{J}:=\tilde{J}\setminus \tilde{I}$.  
Then
$\theta=\hat{I}\cup\hat{J}^\star
\cup M \cup M^\star$ 
and
$\theta^\hash=
=\hat{I}\cup\hat{J}^\star\cup N\cup N^*$.
In other words,  the elements outside of $\hat{I}\cup\hat{J}^\star$
occur in $\theta\cup\theta^\star$ with multiplicity $2$,  those
of $\hat{I}^\star\cup\hat{J}$ with multiplicity $0$,  and the rest
of $\{1,\ldots,2d\}$ with multiplicity $1$.     Now, the right side
of the equation in Proposition~1.8 consists of minors of the form
$(\tilde{J}\cup\bar{\Gamma}, \tilde{I}\cup\bar{\Gamma})$,  where $\bar{\Gamma}$
has the same cardinality as $\Gamma$ and is contained in the complement
of $\tilde{I}\cup\tilde{J}\cup\Gamma$.      Since $\hat{I}$ and $\hat{J}$
remain invariant for all minors appearing in the equation,  it follows
that the multiset $\theta\cup\theta^\hash$ is the same for all terms.

Finally, observe that if $(\theta,\theta^\hash)$ corresponds to an
admissible pair, say $(x,y)$,  then  $\theta\cup\theta^\hash=x\cup y$---see
the proof of Proposition~\ref{pap} above.    Thus when we rewrite a monomial
in the Pl\"ucker coordinates in terms of the minors of \cite{dc},
use the relations of equation~(1.1) and Proposition~1.8 as in the proof
of Theorem~2.4 of~\cite{dc} to express it as a linear combination of
``standard symplectic tableaux'',    and then translate back to get 
a linear combination of standard monomials in our language, the multiset
attached to any monomial on the right side is the same as the multiset
attached to the original monomial. %$\Box$
}.   
%%%%%%%%%End of Footnote %%%%%%%%%%%%%%%%%%%%%%%%%%%%%%%%%%%%%%%%%%%%%
Now let $f_{\varphi_1}\cdots f_{\varphi_s}$ be an arbitrary monomial.  
Consider the expression as a linear combination of standard monomials
for $p_v^h\cdot
p_{\varphi_1}\cdots p_{\varphi_s}$ where $h$ is larger than $2ds$.
It follows from the claim that $p_v$ must occur in every monomial on
the right.   Dividing by $p_v^{h+s}$ we get an expression for
$f_{\varphi_1}\cdots f_{\varphi_s}$ as a linear combination of 
$f_{\ap_1}\cdots f_{\ap_t}$ where $\ap_1\geq\ldots\geq\ap_t$ are $v$-compatible.\eproof

Given an admissible pair $\ap=(x,y)$,  define the {\em $v$-degree} of $\ap$ by
\[ 
\text{$v$-degree}(\ap) := \frac{1}{2}\left( |x\setminus v| + |y\setminus v|
\right) 
\]
An easy calculation using the fact that the $\epsilon$-degrees of $x$ and $y$
are equal gives the following:
\[ 
\text{$v$-degree}(\ap) = |\theta \setminus v| = |\theta^\hash\setminus v|,
\]
where $(\theta,\theta^\hash)$ is the pair associated
to $\ap=(x,y)$ by Proposition~\ref{pap} above.    

The affine patch $G_{d}(V)\cap \mathbb{A}_v$
of the Grassmannian is an affine space whose coordinate ring can
be taken to be the polynomial ring in 
variables of the form $X_{rc}$ where $r$ and $c$ are numbers between $1$ and $2d$
such that $c$ belongs to $v$ and $r$ does not---this is easy to see and 
in any case explained in \cite[\S3]{kr}.  It is readily seen that
the ideal of the closed subvariety 
$\mis\cap\mathbb{A}^v$ of $G_{d,2d}\cap\mathbb{A}_v$
is generated by elements of the form  $X_{rc}-\epsilon X_{c^*r^*}$ where
$\epsilon=-1$ if either $r>d$ and $c^*<d$   or $r<d$ and $c^*>d$,
and $\epsilon=1$ otherwise.     The affine patch $\mis\cap\mathbb{A}^v$
of $\mis$ is thus also an affine space whose  coordinate ring can be
taken to be the polynomial ring in variables of the form $X_\beta$,
$\beta$ in $\roots^v$.
Taking $d=5$ and $v=(3,4,6,9,10)$ for example,   a general element
of $\mis\cap\mathbb{A}_v$ has a basis consisting of column vectors
of a matrix of the following form---having a negative sign
in front of all the variable entries in the first $d$ rows is just
a convenient way of getting the signs right\label{pagematrix}:
\[\left(\begin{array}{ccccc}
-X_{13} & -X_{14} & -X_{16} & -X_{19} & -X_{1,10} \\
-X_{23} & -X_{24} & -X_{26} & -X_{29} & -X_{19} \\
1 & 0 & 0 & 0 & 0 \\
0 & 1 & 0 & 0 & 0 \\
-X_{53} & -X_{54} & -X_{56} & -X_{26} & -X_{16} \\
0 & 0 & 1 & 0 & 0 \\
X_{73} & X_{74} & X_{54} & X_{24} & X_{14} \\
X_{83} & X_{73} & X_{53} & X_{23} & X_{13} \\
0 & 0 & 0 & 1 & 0 \\
0 & 0 & 0 & 0 & 1 \\ \end{array}\right)\]

The expression for $f_\ap:=p_\ap/p_v$ in terms of the $X_\beta$,
$\beta\in\roots^v$, is obtained by taking the determinant of the
 submatrix of a matrix such as above obtained by choosing the
rows given by the entries of $\theta$ where
$\ap\mapsto(\theta,\theta^\hash)$ is the association in
Proposition~\ref{pap}.      Thus it is a homogeneous polynomial
of degree the $v$-degree of $\ap$.    
Since
the ideal of the
Schubert variety $X_w$ in the homogeneous  coordinate ring of
$\mis\subseteq\mathbb{P}(\wedge^d V)$ is generated by 
the $p_\ap$ as $\ap=(x,y)$ varies over all admissible pairs 
such that $w\not\geq x$,\footnote{This is a consequence of Theorem~\ref{tsmt}.
It is easy to see that the $p_\ap$ 
such that $w\not\geq \aptop{\ap}$  vanish on $X_w$.    
Since all standard monomials 
form a basis for the homogeneous  coordinate ring of 
$\mis\subseteq\mathbb{P}(\wedge^d V)$, it follows that 
$w$-dominated standard monomials in admissible pairs span
the quotient ring by the ideal generated by such $p_\ap$.
Since such monomials are linearly independent in the homogeneous
 coordinate ring of $X_w$,    the desired result follows.}
it follows that the ideal of
$\yvw:=X_w\cap\mathbb{A}^v$ in $\mis\cap\mathbb{A}^v$ is 
generated by 
the same $f_\ap$. 
We are interested in the tangent
cone to $X_w$ at $e^v$ (or what is the same, the tangent cone
to $\yvw$ at the origin),   and since $k[\yvw]$ is graded,
its associated graded ring with respect to
the maximal ideal corresponding to the origin
is $k[\yvw]$ itself.   

Proposition~\ref{psmt} tells us that the graded piece of $k[\yvw]$ is
generated as a $k$-vector space by elements of $\smvw$ of degree $m$,
where the degree of a standard monomial $\ap_1\geq\ldots\geq\ap_t$ is
defined to be the sum of the $v$-degrees of $\ap_1,\ldots,\ap_t$.   
To prove Theorem~\ref{tmain} it therefore suffices to prove that
the set $\smvwm$ of $w$-dominated $v$-compatible standard monomials
of degree $m$ is in bijection with $\svwm$.

\section{Further reductions}
\label{sfurtherred}
In the last section we saw how Theorem~\ref{tmain} follows once it is shown
that the combinatorially defined sets $\smvw(m)$ and $\svw(m)$ are in 
bijection with each other.   Now we make some further reductions.  After
these,  it will remain only to show that the combinatorial
bijection established in~\cite[\S4]{kr} has further structure,  
and this will be shown in Section \ref{sproofcompletion}.

An element $v$ of $I(d)$ remains fixed throughout this section.

Let $S^v$ denote the set of monomials in $\roots^v$ and $T^v$ the set of monomials
in $\pos^v$. Let $\smvtopv$ denote the set of $v$-compatible standard monomials 
that are anti-dominated by $v$: a standard monomial $\theta_1\geq \ldots\geq \theta_t$ 
is {\em anti-dominated} by $v$ if $\apbottom{\theta_t}\geq v$. 

Define the {\em domination map} from $T^v$ to $\idn$ by
sending a monomial in $\pos^v$ to the least element
that dominates it\footnote{The poset $I(d)$ has a largest
element, namely $(d+1,\ldots,2d)$, and this clearly dominates
all monomials.  The glb (greatest lower bound) in $I(d,2d)$ of
a set of elements in $I(d)$ belongs to $I(d)$ since $\text{glb}(u,x)^\hash
=\text{glb}(u^\hash,x^\hash)$ for $u$, $x$ in $I(d,2d)$.}.
Define the {\em domination map} from $SM^{v,v}$
to $\idn$ by sending $\theta_1\geq\ldots\geq\theta_t$ to
$\aptop{\theta_1}$. 
Both these maps take, by definition, the value $v$ on the empty monomial.

The desired bijection follows from the following proposition.
\bprop
\label{pbijection}
There is a bijection between $SM^{v,v}$ and $T^v$
that respects domination and degree.
\eprop 

By the above proposition, we 
have, for $w$ in $I(d)$ with $w\geq v$, a bijection $\smvtopvw(m)\cong\tvwm$, 
where $\smvtopvwm$ is the set of $w$-dominated elements of $\smvtopv$
that have degree $m$,   and similarly $\tvwm$ is the set of $w$-dominated
elements of $\tv$ that have degree $m$.

Now let $U^v$ denote the set of monomials in $\roots^v\setminus\pos^v$
(that is,  in pairs $(r,c)$ of $\roots^v$ with $r<c$), and $SM^v_v$ 
the set of $v$-compatible and $v$-dominated standard monomials.  
As explained below,  the ``mirror image'' of the bijection
$SM^{v,v}\cong T^v$   gives a bijection $SM^v_v(m)\cong U^v(m)$,
where $SM^v_v(m)$ and $U^v(m)$ denote respectively the sets of elements
of degree $m$ of $SM^v_v$ and $U^v$.

Putting these bijections together, we get the desired
bijection:
\begin{eqnarray*}
\smvwm &=& 
	\bigcup\limits_{j=0}^m SM^{v,v}_w(j) \times
	SM^v_v(m-j)\\
&\cong&\bigcup\limits_{j=0}^m T^v_w(j) \times
				U^v(m-j)\\
			&=&	 S^v_w(m).
\end{eqnarray*}
Here the first equality is obtained by splitting a $v$-compatible
standard monomial $\theta_1\geq\ldots\geq\theta_t$ into two parts
$\theta_1\geq\ldots\geq\theta_p$ and $\theta_{p+1}\geq\ldots\geq\theta_t$,
where $p$ is the largest integer,  $1\leq p\leq t$,
with $\apbottom{\theta_p}\geq v$. 
The last equality is obtained by writing a
monomial of degree $m$ in $\roots^v$ as a product of two monomials, one
in $\pos^v$ and the other in $\roots^v\setminus\pos^v$,
the sum of  their degrees being $m$.

Let us now briefly explain how to take the ``mirror image''.
Recall that $j^*:=2d-j+1$ for 
an integer $j$, $1\leq j\leq 2d$.
   For $u=(u_1,\ldots,u_d)$ in $\idn$,  define
the dual $u^*$ by $u^*:=(u_d^*,\ldots,u^*_1)$.   This dual map
on $\idn$ is an order reversing involution.  It induces a 
bijection $SM^v_v\cong SM^{v^*,v^*}$ by associating to
$\theta_1\geq\ldots\geq\theta_t$ the element $\theta_t^*
\geq\ldots\geq\theta^*_1$---for an admissible pair
$\theta=(\aptop{\theta},\apbottom{\theta})$, $\theta^*$ denotes of course the
admissible pair $(\apbottom{\theta}^*,
\aptop{\theta}^*)$.   The sum of the $v$-degrees of 
$\theta_1,\ldots,\theta_t$ equals the sum of the $v^*$-degrees of
$\theta_t^*,\ldots,\theta_1^*$,  so that we get a bijection
$SM^v_v(m)\cong SM^{v^*,v^*}(m)$.

For an element $(r,c)$ in $\pos^{v^*}$,  define
the dual to be the element $(c,r)$.   Since $v\in I(d)$, we have
$v^*=[2d]\setminus v$. Thus, it follows that $(c,r)$ belongs to
$\roots^v\setminus\pos^v$.
This induces a degree preserving bijection $T^{v^*}\cong U^v$.
Putting this together with the bijection of the previous paragraph
and the one given by Proposition~\ref{pbijection} 
(for $v^*$ in place of $v$),
we have
\[
	SM^v_v(m)\cong SM^{v^*,v^*}(m)\cong T^{v^*}(m)\cong U^v(m)
\]
The proof of Theorem~\ref{tmain} is thus reduced to the 
the proof of Proposition~\ref{pbijection}.
%%%%%%%%%%%%%%%%%%%%%%%%%%%%%%%%%%%%%%%%%%%%%%%%%%%%%%%%%%%%%%%%%%%%
\subsection{Proof of Proposition~\ref{pbijection}}
The arguments in Section \ref{ssmt} and in the present section thus far
have shown that Theorem~\ref{tmain} follows from
Proposition~\ref{pbijection}.   
To prove the proposition,  we will exploit the work already done in \cite{kr}.
More precisely,  we will deduce the proposition from its earlier version
stated in the paragraphs following Proposition~4.2 of that paper.
As alluded to in the introduction,  the main
ingredient in the deduction is showing that the bijection in the earlier version
respects the involution induced on the combinatorial entities by the
 skew-symmetric form---see Lemma~\ref{lspecial} and Proposition~\ref{pmain}
below for the precise
statements.   We also need Lemma~\ref{ldomination} which roughly speaking is
a symmetry property of domination.
\bremark\label{rkr}
Proposition~\ref{pbijection} and its earlier version are worded
exactly alike.   The difference lies in the meaning attached to the symbols.
To distinguish the two meanings, we use in this section a
tilde sign over these 
symbols to indicate
that the symbol in question has the meaning given to it in \cite{kr}.
However,  in Section \ref{sproofcompletion},  the notation and terminology of that paper
will throughout be in force,  and so the 
$\ \widetilde{ } \ $ will be omitted.
\eremark

Consider $v$ as an element of $I(d,2d)$.  A
{\em standard monomial\/} in $I(d,2d)$ is a totally ordered sequence
$\theta_1\geq\ldots\geq\theta_t$ of elements of $I(d,2d)$.   Such a
monomial is {\em $v$-compatible\/} if each $\theta_j$ is comparable
to $v$ but no $\theta_j$ equals $v$;  it is {\em anti-dominated by $v$\/} if
$\theta_t\geq v$.   

\bnotation\label{nkr} 
In keeping with Remark~\ref{rkr}, we have the following:
\begin{itemize}
\item[ ] 
$\widetilde{\smvtopv}$ denotes the set of $v$-compatible
standard monomials in $I(d,2d)$ anti-dominated by $v$.
\item[ ] $\widetilde{\roots^v}$ denotes the set of ordered pairs $(r,c)$ such that
$r\in [2d]\setminus v$ and $c\in v$.
\item[ ]$\widetilde{\pos^v}$ denotes the subset of $\widetilde{\roots^v}$ consisting
of those $(r,c)$ with $r>c$.
\item[ ] $\widetilde{T^v}$ denotes the set of monomials in $\widetilde{\pos^v}$.
\end{itemize}
\enotation

Consider the natural injection $\smvtopv\to\widetilde{\smvtopv}$: 
$\theta_1\geq\ldots\geq\theta_t$ is mapped to $\aptop{\theta_1}\geq
\apbottom{\theta_1}\geq\ldots\geq\aptop{\theta_t}\geq\apbottom{\theta_t}$, 
but $\apbottom{\theta_t}$ is omitted if it equals~$v$.    This map preserves
domination and doubles degree.   Composing this with the domination 
and degree preserving bijection $\widetilde{\smvtopv} \to \widetilde{T^v}$ of \cite[\S4]{kr},
we get an injection of $\smvtopv$ into $\widetilde{T^v}$.

To describe the image of $\smvtopv$ in $\widetilde{T^v}$,  we introduce
a definition.
For $\alpha=(r,c)$ in $\widetilde{\roots^v}$,
set $\alpha^\hash:=(c^*,r^*)$.
The map $\alpha\mapsto\alpha^\hash$ is an involution on
$\widetilde{\roots^v}$ and on $\widetilde{\pos^v}$
Elements of the form $(r,r^*)$ of $\widetilde{\pos^v}$ 
are referred to as belonging
to the ``diagonal'',  and the set of diagonal elements of $\widetilde{\pos^v}$ 
is denoted $\diagv$.
\bdefinition\label{dspecial} A monomial $\mon$ of $\widetilde{T^v}$ is {\em special\/} if
\begin{enumerate}
\item $\mon=\monhash$  and
\item the multiplicity of any diagonal element in $\mon$ is even.
\end{enumerate}
Equivalently,  $\mon$ is {\em special\/} if there exists $\min$ in
$\widetilde{T^v}$ with $\mon=\min\cup\min^\hash$.
The set of special monomials is denoted $\spec$.
\edefinition

The following is the main technical result of the present paper. 
It is an immediate corollary of Proposition~\ref{pmain} below.
Here we use it to prove
Proposition~\ref{pbijection}.
\blemma\label{lspecial} 
The image of $\smvtopv$ in $\widetilde{T^v}$ is the set $\spec$ of
special monomials. \elemma
In other words,    there is a domination preserving and degree doubling
bijective map from $\smvtopv$ to $\spec$.  

On the other hand,
there is a domination preserving and degree halving bijective map from
$\spec$ to $T^v$.   Given $\mon$ in $\spec$, to get the corresponding
element of $T^v$, 
replace those $(r,c)$ of $\mon$ with $r>c^*$ by $(c^*,r^*)$ and then
take the (positive) square root.     This clearly halves the degrees.
The map in the other direction is obvious:  given a monomial of $T^v$, 
replace each $(r,c)$ by $(r,c)(c^*,r^*)$---in other words,  a monomial
$\frak P$ of $T^v$ is mapped to $\frak P \cup {\frak P}^\hash$.

To see that this bijection from $\spec$ to $T^v$ preserves domination,
let $\frak T$ be a special monomial and $\frak P$ the corresponding
monomial in $T^v$.   Let $y$ and $z$ be the images respectively of 
$\frak T$ and $\frak P$ under the domination maps.  Then $y$ is the
first  coordinate of $\pi(\frak T)$,  where $\pi$ is the map defined
in \cite[\S4]{kr}---see Proposition~4.1~(4) of \cite{kr}.  Since
$\frak T ={\frak T}^\hash$,  it follows from Proposition~\ref{ppihashrespect}
below that $y=y^\hash$.   This is equivalent to saying that $y$
belongs to $I(d)$.

That $y$ dominates $\frak P$ is clear: any $v$-chain in $\frak P$ is also
a $v$-chain in $\frak T$.     Since $y$ is in $I(d)$,  it follows from
the definition of the domination map on $T^v$ that $y\geq z$.   That
$z\geq y$ follows from Lemma~\ref{ldomination} below: $z^\hash=z$,
$\frak P={\frak P}^{\rm{up}}$, and $z$ dominates $\frak P$, so that
$z$ dominates $\frak P\cup{\frak P}^\hash=\frak T$.

Composing the two bijective maps $\smvtopv\to\spec$ and $\spec\to T^v$ above
gives us a domination and degree preserving bijection
$\smvtopv\to T^v$.    Proposition~\ref{pbijection} is thus proved.

\section{Completion of proof}\label{sproofcompletion}
The purpose of this section is to establish the statements used in the proof
of Proposition~\ref{pbijection}.   The terminology and notation of \cite{kr}
will throughout be in force.

An element $v$ of $I(d,n)$ remains fixed throughout. 
The symbol $w$ will denote an arbitrary element of $I(d,n)$ satisfying $w\geq v$, 
$\mon$ an arbitrary monomial in $\pos^v$, and
$\monw$ the unique subset of $\pos^v$
defined by Proposition~\ref{poldmonw} below.   

Let $(r,c)$ and $(R,C)$ be elements of $\pos^v$.   They are {\em comparable}
if either they are equal or $(R,C)>(r,c)$ (which means that $R>r$ and $C<c$)
or $(r,c)>(R,C)$.    We say that $(R,C)$ {\em dominates} $(r,c)$ if
$R\geq r$ and $C\leq c$.

\subsection{Recall}
\label{ssrecall}\label{ssmonw}
\bprop
\label{poldmonw}
To each element $w$ of $I(d,n)$ satisfying $w\geq v$ there is associated a unique
subset $\mon_w=\{(r_1,c_1),\ldots,(r_p,c_p)\}$ of $\pos^v$ satisfying
conditions A--C below.   
\begin{enumerate}
\item[A.]
 $r_i\neq r_j$ and $c_i\neq c_j$ for $i\neq j$.
\item[B.] If $r_i<r_j$ then either $c_j<c_i$ or $r_i<c_j$. 
\item[C.]  $w =\left(v\setminus\{c_1,\ldots,c_p\}\right)\cup\{r_1,\ldots,r_p\}$.
\end{enumerate}
Furthermore, $\mon_w$ also satisfies conditions D and E below.
\begin{enumerate}
\item[D.] The $v$-degree of $w$ equals the cardinality $p$ of $\mon_w$.
\item[E.] $w$ is the smallest element of $I(d,n)$ to dominate $\mon_w$.
\end{enumerate}
\eprop
\version{   %%%%%%%%%%%%%%%%%%%  Public version   %%%%%%%%%%%%%%%%%%%
\bproof  The existence and uniqueness of $\mon_w$ is the substance
of Proposition~4.3 of \cite{kr}.   That D holds is obvious.
For E,  see Lemma~5.5 of \cite{kr}.\eproof}
{ %%%%%%%%%%%%%%%%%%%%   private version   %%%%%%%%%%%%%%%%%%%%%%%%%
}
\bremark
If $\mon$ is a subset of $\pos^v$ that satisfies condition A of
Proposition~\ref{poldmonw},   the equation in item C can be taken to be
the definition of an element $w$ of $I(d,n)$.    If $\mon$ satisfies $B$
also,    then $\mon=\mon_w$ for this $w$.    Let $x$ be an element of $I(d,n)$
and $\mon$ a subset of $\mon_x$.    Then $\mon$ evidently satisfies A and B.
Therefore $\mon=\mon_w$ for some $w$.    It follows from E that $x\geq w$.
\eremark
%%%%%%%%%%%%%%%%%%%%%%%%%%%%%%%%%%%%%%%%%%%%%%%%%%%%%%%%%%%%%%%%%%%%%%%%%%%
\blemma\label{lblockdisjoint}
Let $\block$ and $\blc$ be distinct blocks of a monomial $\mon$.
Then $\block\cup\block'$ and $\blc\cup\blc'$ are disjoint.
\elemma
\bproof
For $\block$ and $\blc$ both belonging to the same partition $\mon_j$
of $\mon$,   this follows from the proof of Lemma~4.10 of \cite{kr}---see
the last line of that proof.     Now suppose that $\block\subseteq\mon_i$
and $\blc\subseteq\mon_j$,   and that $(r,c)$ belongs to both 
$\block\cup\block'$ and $\blc\cup\blc'$.    Assume without loss of
generality that $i<j$.   From the definition of $\blc'$,  it follows that there
exists an element $(r',c)$ in $\blc$ with $r'\geq r$.   From the definition
of the partitions $\mon_k$,   it follows that there exists $(a,b)$ in 
$\mon_i$ such that $(a,b)>(r',c)$.   Now $(a,b)>(r,c)$ and both of them
belong to $\mon_i\cup\mon_i'$,   a contradiction to Lemma~4.10 of \cite{kr}.\eproof
\bcor\label{clblockdisjoint}
Given $(r,c)$ in $\mon\cup\mon'$,   there exists a unique block $\block$ of
$\mon$ such that $(r,c)$ belongs to $\block\cup\block'$.$\Box$ \ecor
%%%%%%%%%%%%%%%%%%%%%%%%%%%%%%%%%%%%%%%%%%%%%%%%%%%%%%%%%%%%%%%%%%%%%%%%%%%
\subsection{Statement of the main proposition}
\label{ssstatement}\label{ssstatements}
From now on we will assume that $n=2d$ and that $v=v^\invol$ (or equivalently, 
that $v$ belongs to $I(d)$).
The map $\pi$ of \cite[\S4.2]{kr} is defined only on nonempty monomials.
We extend the definition by setting $\pi(\mon)=(v,\mon)$ if $\mon$ is
empty.   Similarly we extend the definition of the map $\phi$ of
\cite[\S4.4]{kr} by setting
$\phi(v,\min)=\min$ if $\min$ is empty.

The following proposition has as an immediate corollary
Lemma~\ref{lspecial} which was the main ingredient in the proof of 
Proposition~\ref{pbijection}.    Its proof will be given in
Subsections \ref{ssproofmainone} and \ref{ssproofmaintwo}
after some preliminaries in Subsection \ref{sshashrespect}.

\bprop\label{pmain}
Let $\mon$ be a special monomial.   Set 
$$\pi(\mon)= (t,\mon') \quad \mbox{ and } \quad \pi(\mon')=(u,\min).$$ 
Then
\begin{enumerate}
\item
$t\geq u\geq v$,
$t^\invol =t $,  $u^\invol =u$,
and the $(1,\ldots,d)$-degrees of $t$ and $u$ are equal.
\item
$\min$ is special and $u$ dominates $\min$.
\end{enumerate}
Conversely,   if $t$ and $u$ are elements of $I(d,2d)$ and $\min$ a
monomial such that the two conditions above are
satisfied,   then $\phi(t,\phi(u,\min))$ is special.  
\eprop
%%%%%%%%%%%%%%%%%%%%%%%%%%%%%%%%%%%%%%%%%%%%%%%%%%%%%%%%%%%%%%%%%%%%%%%%%%%
\subsection{Behavior under the operation $\invol$}\label{sshash}\label{sshashrespect}
In this subsection we investigate the
behavior of the combinatorial constructions of \cite[\S4]{kr} under the
operation $\invol$.

If $\beta$ is $j$-deep in $\mon$, then $\beta^\invol$ is $j$-deep in $\mon^\invol$.
Thus the partitions $\mon_j$ of $\mon$ are respected by the hash operation:
$(\mon_j)^\invol = (\mon^\invol)_j$.    It is also easy,
given the definitions,
to verify the following:
\begin{itemize}
\item If $\block$ is a block of $\mon_j$,  then $\block^\invol$ is a block of
	$\mon_j^\invol$.
\item $w(\block^\invol) = w(\block)^\invol$.
\item $(\block^\invol)' = (\block')^\invol$.
\end{itemize}
These observations amount to a proof of the following proposition.

\bprop\label{ppihash}\label{ppihashrespect}
The map $\pi$ respects $\invol$. More precisely, if $\pi(\mon)=(w,\mon')$,
then $\pi(\mon^\invol)= (w^\invol,{\mon'}^\invol)$.
\eprop

Let us now show that the map $\phi$ also respects $\invol$.  For this we need

\bprop\label{pmonwhashrespect}\label{pmonwhash}  
The association $w\mapsto\monw$ respects $\invol$, that is,  
$(\monw)^\invol = \mon_{w^\invol}$.   In particular,  $w^\invol=w$
if and only if $(\monw)^\invol=\monw$.
\eprop

\bproof   
The result follows from Proposition~\ref{poldmonw}:   $(\monw)^\invol$
satisfies conditions A, B and C of that proposition with 
$w$ replaced by $w^\invol$ in condition C.
\eproof

Let $\min$ be a monomial in $\pos^v$ and $w$ an element of $I(d,2d)$
such that $w\geq v$ and $w$ dominates $\min$.
The following are easily verified from the definitions:
\begin{itemize}
\item $(\monw^j)^\invol = \mon^j_{w^\invol}$ 
and so $(w^j)^\invol =(w^\invol)^j$.
\item $(\min_j^w)^\invol = (\min^\invol)^{w^\invol}_j$.
\item $(\piece_\beta)^\invol$ is the piece of $\min^\invol$ corresponding
to $\beta^\invol$,  for $\beta$ an element of $\monw$.
\item
$\left((\min^w_j)^\star\right)^\invol = \left((\min^\invol)_j^{w^\invol}\right)^\star$.
\end{itemize}
These observations amount to a proof of the following proposition.

\bprop\label{pphihashrespect}
The map $\phi$ respects $\invol$:   more precisely,  if 
$\phi(w,\min)= \mon$  then $\phi(w^\invol,\min^\invol)= \mon^\invol$.\eprop

\bprop\label{pdegree}
Suppose that $w=w^\invol$ (or equivalently, by (3), $\mon_w=\mon_w^\invol$).  Then
for an element $(r,c)$ of $\mon_w$, either {\rm (i)}~$c<r\leq d$, or
{\rm (ii)}~$d<c<r$, or {\rm (iii)}~$r^\star=c\leq d <r$.   In particular, 
the $\epsilon$-degree of $w$ equals the number of elements of $\mon_w$
belonging to the diagonal.
\eprop
\bproof 
Suppose that $c\leq d$ and $r>d$.    Since $\monw=\monw^\hash$,
$(c^*,r^*)$ belongs to $\monw$.   If $r<c^*$,  then $r^*>c$,   and the presence
of both $(r,c)$ and $(c^*,r^*)$ in $\monw$ is a violation of condition~E
of Proposition~\ref{poldmonw},  a contradiction.   If $c^*<r$,  then
$c>r^*$,  and again there is a similar contradiction.   Thus $c=r^*$.
\eproof

\bprop\label{pbbhash}\label{pmult}
Let $\block$ be a block of a monomial $\mon$ satisfying $\mon=\mon^\invol$.
\begin{enumerate}\item[(A)] The following are equivalent:
\begin{enumerate}
\item[(1)]	$\block=\bhash$.
\item[(2)]	$w(\block)$ lies in the diagonal.
\item[(3)]	$\block\cup\block'$ meets the diagonal.
\end{enumerate}
\item[(B)]
Suppose that the conditions in (A) above are met. 
Then there is a unique element of the diagonal in $\block\cup\block'$ and the
multiplicities of that element  
in $\block$ and $\block'$ differ
by $1$.
\end{enumerate}
\eprop
\bproof 
{\rm (A):\ }  As noted in the proof of Proposition~\ref{ppihashrespect},
$\block^\invol$ is a block of $\mon^\invol=\mon$,  $w(\block^\invol)=w(\block)^\invol$,
and $(\block^\invol)'=(\block')^\invol$.
\begin{description}
\item[$(1)\Rightarrow(2)$] $w(\block)^\invol=w(\block^\invol)=w(\block)$,  which
means that $w(\block)$ is on the diagonal.
\item[$(2)\Rightarrow(1)$] We have $w(\block^\invol)=w(\block)^\invol=w(\block)$.
If $\block$ and $\block^\invol$ are distinct blocks,   then $w(\block)$ and
$w(\block^\invol)$ cannot share a row or column
index---Corollary~4.13 of \cite{kr}---let
alone being equal.
\item[$(3)\Rightarrow(1)$]
It follows from the hypothesis that $\block\cup\block'$ and $\block^\invol\cup
(\block^\invol)'$ meet.   By Corollary~\ref{clblockdisjoint},  $\block=\block^\hash$.
\item[$(1)\Rightarrow(3)$]
Since $\block=\block^\invol$, there is an element $(r,s)$ of
$\block$ with $r\leq s^*$.   Let $r$ be the maximal row index of such an element
of $\block$.
A portion of the arrangement of the elements of
$\block$ in ascending order of row and column indices looks like this:
\[
	\ldots,(a,b), \underbrace{(r,r^\star),\ldots,(r,r^\star)}_{\text{$m$ times}},
	(b^\star,a^\star),\ldots	
\]
where either $a<r$ or $b<r^\star$.
If $m>0$ then $\block$ meets the diagonal and we are done.
If $m=0$,   then $a=r$, and the portion of $\block'$ corresponding to the one
above of $\block$ looks like this:
\[
	\ldots,(?,b),(a,a^\star),(b^\star,?),\ldots	\]
so that $(a,a^\star)=(r,r^\star)$ belongs to $\block'$.
\end{description}
{\rm(B):\ } Since elements of the diagonal are comparable but no two distinct elements
of $\block\cup\block'$ are (the latter statement is immediate from the
definitions; or see Lemma~4.10 of \cite{kr}),   there is at most one diagonal
element in $\block\cup\block'$.   Since there is at least one such element
by hypothesis,  there is a unique such element.

Let $(r,r^\star)$ be the unique diagonal element $\block\cup\block'$.  Proceed
as in the proof of $(1)\Rightarrow(3)$ above.   If $m=0$, then  
$(a,a^*)=(r,r^*)$ belongs to $\block'$.   Further $b<r^*$ and $b^*>r$,  so that
the multiplicity of $(r,r^\star)$ in $\block'$ is $1$ and we are done.

If $m>0$,  
the portion of $\block'$ corresponding to the one
above of $\block$ looks like this:
\[
\ldots,(?,b),(a,r^\star),\underbrace{(r,r^\star),\ldots, (r,r^\star)}_{\text{$m-1$ times}},
	(r,a^\star),(b^\star,?),\ldots	\]
The multiplicity of $(r,r^\star)$ in $\block'$ is therefore either $m+1$ or $m-1$
depending upon whether or not $a$ equals $r$.
\eproof

\bcor\label{cpmult}
Let\/ $\mon=\mon^\invol$ and $(r,r^\star)$ be an element of the
diagonal
belonging to $\mon\cup\mon'$.   Then the multiplicities of $(r,r^\star)$ in
$\mon$ and $\mon'$ differ by $1$.
\ecor

\bproof
By Corollary~\ref{clblockdisjoint}, there is a unique block
$\block$ of $\mon$ such that $(r,r^\star)$ belongs to $\block\cup\block'$.  The
multiplicities of $(r,r^\star)$ in $\mon$ and $\mon'$ equal respectively those
in $\block$ and $\block'$,   and 
these differ by $1$ by the proposition above.
\eproof

\bdefinition
\label{dupdown} 
For $\alpha=(r,c)$ in $\pos^v$, set
\[ \alpha^{\rm up} := \left\{ \begin{array}{cl}
	\alpha & \mbox{ if $r\leq c^\star$}\\
	\alpha^\invol =(c^\star,r^\star) & \mbox{ if $r>c^*$}
	\end{array}\right.
	\quad \ \ 
 \alpha^{\rm down} := \left\{ \begin{array}{cl}
	\alpha & \mbox{ if $r\geq c^\star$}\\
	\alpha^\invol =(c^\star,r^\star) & \mbox{ if $r<c^*$}
	\end{array}\right. 	
\]
\edefinition

\bprop
\label{pdiagbelong}
Let $\alpha$ be an element of $\pos^v$ and $w$ an element of $I(d,2d)$
with $w\geq v$.   If both 
$\alpha^{\rm up}$ and $\alpha^{\rm down}$ belong to $\monw$
and both of them dominate an element $\beta$ of $\pos^v$,
then $\alpha^{\rm up}=\alpha^{\rm down}$.    In particular, if $w=w^\invol$
and $\alpha$ in $\monw$ dominates an element on the diagonal,
then $\alpha$ belongs to the diagonal.
\eprop

\bproof  
Write $\alpha^{\rm up}=(R,C)$ and $\beta=(r,c)$.  Then
$C\leq R^*\leq c<r\leq R\leq C^*$---here $C\leq R^*$ and $R\leq C^*$ because
$\alpha^{\rm up}=(R,C)$,   $R^*\leq c$ because $\alpha^{\rm down}=(C^*,R^*)$
dominates $\beta$,  $c<r$ because $\beta$ belongs to $\pos^v$, and 
$r\leq R$ because $\alpha^{\rm up}$ dominates $\beta$.
Unless $R=C^*$
the presence of both $\alpha^{\rm up}$ and $\alpha^{\rm down}$ in $\monw$
leads to a violation of condition~B of Proposition~\ref{poldmonw} by $\monw$,
a contradiction.
\eproof

\bcor\label{cpdiagbelong}
If $\alpha_1>\ldots>\alpha_p$ is the $v$-chain of diagonal
elements in $\monw$ for an element $w=w^\invol$,  
then $\alpha_j$ is $j$-deep but
not $j+1$-deep in $\monw$.$\Box$
\ecor

\blemma\label{ldompreserve}
\label{ldomination}
If\/ $w^\hash=w$ and $w$ dominates a monomial $\mon$ with $\mon=\mon^\up$
(that is, $\alpha =\alpha^\up$ for every $\alpha\in\mon$),  then $w$
also dominates $\mon\cup\mon^\hash$. 
\elemma

\bproof   
Let $\beta_1>\ldots>\beta_p$  be a $v$-chain in
$\mon\cup\mon^\hash$.     We need to show that this chain is dominated by
a $v$-chain in $\mon_w$.    
We may assume without loss of generality
that $\beta_p^\up=\beta_p$.   If $\beta_j^\up=\beta_j$  for all $j$,
then the desired result follows immediately from the hypothesis.  
Otherwise, let $k\in [p]$ be the largest such that $\beta_k^\down=\beta_k$.
Choose $v$-chain $\alpha_1>\ldots>\alpha_p$ in $\mon_w$ such that
it dominates $\beta_1^\up>\ldots>\beta_p^\up$.   
Replacing $\alpha_1>\ldots>\alpha_p$ by $\alpha_1^{\rm up}>\ldots>\alpha_p^{\rm up}$,
we may assume that $\alpha_j^\up=\alpha_j$.
It follows from Proposition~\ref{pdiagbelong} that $\alpha_j$
lives on the diagonal for $1\leq j\leq k$:   $\alpha_k^{\rm down}>
\beta_k^{\rm down}=\beta_k>\beta_{k+1}$ and $\alpha_k^{\rm up}>
\beta_k^{\rm up}>\beta_{k+1}^{\rm up}=\beta_{k+1}$  implies $\alpha_k$
belongs to the diagonal;   and $\alpha_{k-1}>\alpha_k$ implies 
that $\alpha_{k-1}$ belongs to the
diagonal etc.
Thus $\alpha_j$ dominates $\beta_j$ for all $j\in [p]$. %, $1\leq j\leq p$.
\eproof

\subsection{Proof of the first half of  Proposition~\ref{pmain}}
\label{ssproofone}\label{ssproofmainone}

By Proposition~\ref{ppihashrespect} above, $\pi$ respects the $\invol$
operation.
Hence it follows that 
$t^\invol=t$, $(\mon')^\invol = \mon'$, $u^\invol=u$ and $(\min)^\invol = \min$.   
Also, it follows from 1, 3~and~4 of~Proposition~4.1 of \cite{kr} that $t\geq u\geq v$, 
and that $u$ dominates $\min$.  
Thus it remains only to prove the following:
\begin{enumerate}
\item[(a)] $\mon_t$ meets the diagonal in as many points as $\mon_u$---this
will imply,  by Proposition~\ref{pdegree},  that the $\epsilon$-degrees
of $t$ and $u$ are equal and hence that $(t,u)$ is an admissible pair.
\item[(b)] $\min$ is special.
\end{enumerate}

We first prove (a).  
If $(r,r^\star)$ and $(s,s^\star)$ are elements of the diagonal
$\diagv$,   then either $(r,r^\star)>(s,s^\star)$ or the other way around.  
The elements 
of $\mon_u\cap\diagv$ therefore form a $v$-chain.  Let
$\beta_1>\ldots>\beta_p$ be all the elements of $\mon_u\cap\diagv$.  
By Proposition~\ref{poldmonw}~(E),  $t$ dominates $\mon_u$.
By Lemma~4.5 of \cite{kr},  there exists a $v$-chain $\alpha_1
>\ldots > \alpha_p$ in $\mon_t$ that dominates $\beta_1>\ldots>\beta_p$.
By Proposition~\ref{pdiagbelong},   the $\alpha_j$ belong to $\diagv$.

To complete the proof of (a), it remains to be seen that
the cardinality of $\mon_u\cap\diagv$ is not less than
that of $\mon_t\cap\diagv$.    We will show in fact that $\mon'\cap\diagv$
has at least as many distinct elements as the cardinality of
$\mon_t\cap\diagv$.      The blocks containing these elements will then
clearly be distinct:  elements of the diagonal are comparable but no two distinct
elements of a block are.    Each such block $\blockc$ satisfies
$\blockc^\invol=\blockc$, for blocks are disjoint and $\blockc$ meets
the diagonal.    If $\blockc$ and $\blockc_1$ are distinct blocks,
then $w(\blockc)$ and $w(\blockc_1)$ are distinct---Corollary~4.13 of
\cite{kr}---so there are distinct contributions to $\mon_u$ from these
blocks.   Finally,   these contributions all lie on the diagonal
by Proposition~\ref{pbbhash}.

Let $\beta$ be an element of $\mon_t\cap\diagv$.
Let $\block$ be the block of $\mon$ with $w(\block)=\beta$.
By Proposition~\ref{pbbhash},  $\block\cup\block'$ meets
$\diagv$.    If $\block$ meets $\diagv$ at say $\beta'$, then,
since $\mon$ satisfies condition~2 of Definition~\ref{dspecial},
   it follows from the definition of $\block'$ that
$\block'$ also contains $\beta'$.    Thus we get an element 
$\beta'$ in $\mon'\cap\diagv$ for every $\beta$ in $\mon_t\cap\diagv$.
The association $\beta\mapsto\beta'$ is injective, for
the association $\beta\mapsto\block$ is injective---see the reason
given in the last paragraph---and so is $\block\mapsto\beta'$ by
Lemma~\ref{lblockdisjoint}.    This finishes the proof of (a).

For use in the proof of (b), let us record a corollary of the proof
of (a):
\bcor\label{cpmainone}   
A block $\block$ of $\mon'$ meets the diagonal
if and only if $\block=\block^\invol$.\ecor

\bproof   
If $\block$ meets the
diagonal,  then $\block^\invol=\block$.   The number of 
$\block$ such that $\block^\invol=\block$ is the cardinality of
$\mon_u\cap\diagv$.    The number of $\block$ that meet the diagonal
is the cardinality of $\mon'\cap\diagv$.   But we have seen that the
cardinalities of $\mon_u\cap\diagv$, $\mon'\cap\diagv$, and $\mon_t\cap
\diagv$ are equal.
\eproof

We now turn to the proof of (b).    Suppose that $(r,r^\star)$ belongs
to $\min$.   Let $\block$ be the block of $\mon'$ such that
$(r,r^\star)$ belongs to $\block'$.   By Proposition~\ref{pbbhash} and
the corollary above,  $\block$ meets the diagonal.   Since elements of
$\block\cup\block'$ are incomparable---see Lemma~4.10 of \cite{kr}---it 
follows that $(r,r^\star)$ belongs to $\block$  and so to $\mon'$.
By Corollary~\ref{cpmult},   the multiplicities of $(r,r^\star)$ in 
$\mon$ and $\mon''=\min$ differ by $2$ or $0$ or $-2$. %\hfill $\Box$

\subsection{Proof of the converse in  
Proposition~\ref{pmain}}\label{ssprooftwo}\label{ssproofmaintwo}
Set $\monm=\phi(u,\min)$ and $\monn=\phi(t,\monm)$.   Since $\phi$
respects $\hash$---see Proposition~\ref{pphihashrespect}---it follows
that $\monm^\invol=\monm$  and $\monn^\invol=\monn$.   The only thing
that takes some work is to show that $\monn$ satisfies condition~2 of
Definition~\ref{dspecial}.

By Proposition~\ref{pbbhash}, each diagonal element of
$\monn$ belongs to a block $\block$ of $\monn$ with $w(\block)$ on the
diagonal.    It therefore suffices to show that the multiplicity of
a diagonal element in any such block of $\monn$ is even (possibly $0$).

Let $\alpha_1>\ldots>\alpha_p$ and $\beta_1>\ldots>\beta_p$ be respectively
the elements $\mon_t\cap\diagv$ and $\mon_u\cap\diagv$.
Let $\piece_1$, \ldots, $\piece_p$ be the pieces of $\min$ with respect
to $u$ corresponding respectively to $\beta_1$, \ldots, $\beta_p$.  Then
$\piece_j^\star=\block$ is a block of $\monm$ with
$w(\piece_j^\star)=\block_j$ and $\piece_j=\block'$---see
the claim at the beginning of \S4.5 of \cite{kr}. 

Since $\min$ is special,  it follows from Proposition~\ref{pbbhash}
that $\piece_j^\star$ meets the diagonal with (positive) odd multiplicity.
Denoting the diagonal element of $\piece_j^\star$ by $\delta_j$,
we have $\delta_1>\ldots>\delta_p$---see Lemma~4.21 of \cite{kr}.  
By Proposition~\ref{pdiagbelong},
the piece of $\monm$ with respect to $t$ to which $\delta_j$ belongs
can only correspond to one of the $\alpha_i$.

No two distinct elements of a piece are comparable,  so two different $\delta$ cannot
belong to the piece corresponding to the same $\alpha$.  Further,  it
follows from the the definition of the way we break a monomial into pieces
that if $i<j$ and $\delta_i$, $\delta_j$ belong respectively to the pieces
corresponding to $\alpha_k$ and $\alpha_l$,   then $k<l$.   Thus $\delta_j$
belongs to the piece corresponding to $\alpha_j$.    

Let $\block_1$, \ldots, $\block_p$ be the blocks of $\monn$ with
$w(\block_j)=\alpha_j$.      Then $\block_j'$ is the piece of $\monm$ with
respect to $t$ corresponding to $\alpha_j$---see the claim at the beginning of
\cite[\S4.5]{kr}.     By
Proposition~\ref{pbbhash},  the multiplicity of 
the diagonal element in $\block_j$ is even (possibly $0$). \hfill $\qed$

\section{Interpretations}\label{sinterpretations}
Fix elements $v,w$ in $I(d)$ with $v\leq w$. 
It follows from Corollary~\ref{cmain} that the
multiplicity of the Schubert variety $X_w$ in
$\mis$ at the point $e^v$ can be
interpreted as the cardinality of a certain set
of nonintersecting lattice paths.    We 
illustrate this by means of two examples.   The justification 
for the interpretation is the same word-for-word as for the usual
Grassmannian~\cite[\S5]{kr} and so will be omitted. The interpretation
leads immediately to 
a formula for the multiplicity involving binomial determinants.
Specializing to the case $v=\epsilon = (1,\ldots,d)$,  we recover
Conca's result~\cite[Theorem~3.6]{conca}.

\bexample\label{eone}
Let $d=23$. Consider
$$
v=(1,2,3,4,5,11,12,13,14,19,20,22,
23,26,29,30,31,32,37,38,39,40,41)
$$
and
$$
w=(4,5,9,10,14,17,18,21,23,25,27,28,
31,32,34,35,36,39,40,41,44,45,46) 
$$
so that 
\[\begin{array}{rcl}
\mon_w & = &\{(9,3), (10,2), (17,13), (18,12), (21,20),  (25,22),\\
	& & \; (27,26), (28,19), (34,30), (35,29), (36,11), (44,38), (45,37), 
(46,1)  \}. 
\end{array}\]

\begin{figure}\label{fexone}
\begin{center}
\mbox{\epsfig{file=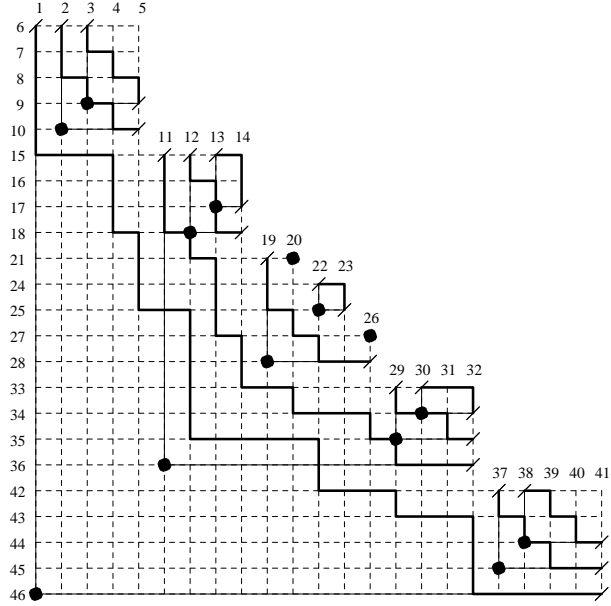,height=8cm,width=8cm}}
\caption{\label{lfexone}A symmetric tuple of nonintersecting lattice paths as in 
	Example~\ref{eone}}
\end{center}
\end{figure}

We think of $v$ and $w$ as elements of $I(d,2d)$ and switch to the
notation of~\cite{kr} as we have done in \S\ref{sproofcompletion}.
The grid depicting the points of $\pos^v$ is shown
in Figure~\ref{fexone}.  The solid dots represent the
points of $\mon_w$, this being the monomial associated to $w$ by
Proposition~\ref{poldmonw}.   From each point $\beta$
of $\mon_w$ we
draw a vertical line and a horizontal line.
Let $\bstart$ and $\bfinish$ denote respectively the points
where the vertical line and the horizontal line meet
the boundary.     For example, $\bstart=(15,11)$ and
$\bfinish=(36,32)$ for $\beta=(36,11)$;   for $\beta=(21,20)$,
$\bstart=\bfinish=\beta$.

A {\em lattice path} between a pair of such points $\bstart$
and $\bfinish$ is a sequence $\alpha_1,\ldots,\alpha_q$
of elements of $\pos^v$ with $\alpha_1=\bstart$ and
$\alpha_q=\bfinish$  
such that for $1\leq j\leq q-1$, if we let $\alpha_j=(r,c)$, then 
$\alpha_{j+1}$ is either 
%$(r+1,c)$ or $(r,c+1)$.    %<----mistake; corrected below
$(R,c)$ or $(r,C)$ where $R$ is the least element of $[2d]\setminus v$
that is bigger than $r$ and $C$ the least element of $v$ that is bigger
than $c$.
Note that if  $\bstart=(r,c)$ and
$\bfinish=(R,C)$,   then  %$q=(R-r)+(C-c)+1$.  %<----mistake; corrected below
$q$ equals
\[ |([2d]\setminus v)\cap\{r,r+1,\ldots,R\}| +
 | v\cap\{c,c+1,\ldots,C\}| - 1. \]

Write $\mon_w=\{\beta_1,\ldots,\beta_p\}$. 
By Proposition~\ref{pmonwhash}, $\mon_w^\hash=\mon_{w^\hash}=\mon_w$---that
is,  $\mon_w$ is symmetric with respect to the anti-diagonal.
Consider
the set of all $p$-tuples of paths $(\path_1,\ldots,\path_p)$,
where $\path_j$ is a lattice path between $\bstartj$
and $\bfinishj$ such that 
\begin{itemize}
\item  No two $\path_j$ intersect.
\item  If $\beta_j=\beta_k^\hash$,  then $\Lambda_j=\Lambda_k^\hash$.
\end{itemize}
A particular such symmetric $p$-tuple is shown in Figure~\ref{fexone}.
The number of such symmetric $p$-tuples is the multiplicity of 
$X_w$ at the point $e^v$.
\eexample

\bexample\label{etwo}
Let us draw, in a simple case,
the pictures of all possible symmetric tuples of nonintersecting
lattice paths as defined in the above example. 
Let 
$$
d=5, \quad v=(1,2,3,6,7) \quad \mbox{ and } \quad 
w=(3,5,7,9,10)  
$$
so that $\mon_w=\{(5,2), (9,6), (10,1)\}$.
Figure~\ref{fextwo} shows all
the symmetric $3$-tuples of 
nonintersecting lattice paths.  There are $10$ of them
and thus the multiplicity in this case is $10$.
\eexample

\subsection{The Gr\"obner basis interpretation}\label{sgroebner}\label{sgrobner}
Here we interpret
Theorem~\ref{tmain} 
in terms of Gr\"obner basis---see Proposition~\ref{pgrob}
below.     The special case of the proposition when $v$ is
the identity coset (that is, $v=(1,\ldots,d)$)
has been obtained by Conca~\cite{conca} by different methods.

Fix elements $v,w$ in $I(d)$ with $v\leq w$. 
Recall from \S\ref{ssmt} that the ideal
$$
\left(f_\ap=p_\ap/p_v\st
\text{ $\ap$ admissible pair, $\aptop{\ap}\nleq w$}\right)
$$ 
in the polynomial ring $P:=k[X_\beta\st\beta\in\roots^v]$ defines the 
tangent cone to the Schubert variety $X_w$ at the point $e^v$.   
We will identify a subset of these generators~$f_\ap$ as being Gr\"obner
basis with respect to certain term orders for the ideal of the tangent cone.
Observe that the $f_\ap$ are homogeneous polynomials in the
variables $X_\beta$.

\begin{figure}
\begin{center}
\mbox{\epsfig{file=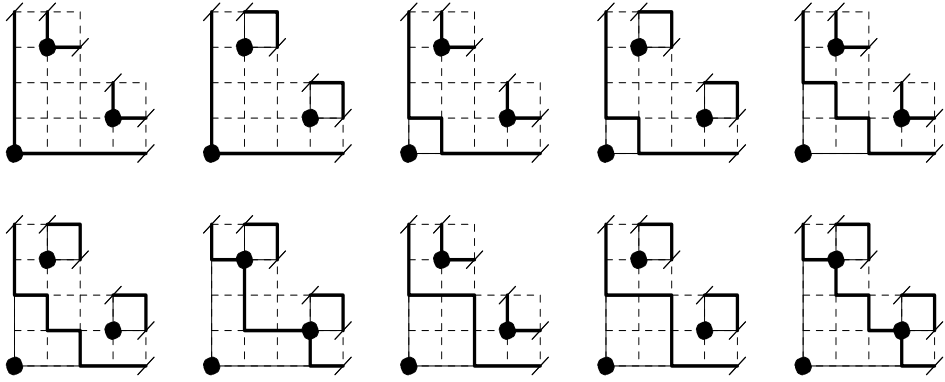,height=5cm,width=12.5cm}}
\caption{\label{fextwo}All the symmetric tuples of nonintersecting 
lattice paths as in Example~\ref{etwo}}
\end{center}
\end{figure}

Let $\ap=(t,u)$ be an admissible pair and $(\theta,\theta^\hash)$
the pair associated to $\ap$ as in Proposition~\ref{pap}.  
Let us call $\ap$ ``good'' (the scope of this
terminology is intended to be limited to this section) if 
\begin{itemize}
\item $v\leq u$ and $t\nleq w$;
\item $\mon_\theta=\mon_\theta^\text{up}$,  that is,
$r\leq c^*$ for $(r,c)$ in $\mon_\theta$,  and the elements of 
$\mon_\theta$ form a $v$-chain---here $\mon_\theta$ is the monomial
attached to $\theta$ as in Proposition~\ref{poldmonw}.\footnote{
While it is true that in
\S\ref{sproofcompletion} we are using the notation of \cite{kr}
where the symbols have a different meaning from what they do in \S\ref{stheorem}--\ref{sfurtherred},
the condition that $r\leq c^*$ for $(r,c)$ in $\mon_\theta$ gives
us some license for sloppiness.}
\end{itemize}

Suppose that $\ap$ is good.   Then $\mon_\theta$ occurs up to sign
as a term in the expression for $f_\ap=f_\theta$ as the determinant
of a submatrix of the matrix of the form on page~\pageref{pagematrix}
attached to $v$---here we are abusing notation and not distinguishing
between a monomial in $\roots^v$ and the corresponding monomial
in the variables $X_\beta$, $\beta$ in $\roots^v$.

\bprop\label{pgrob}\label{pgroebner}
Fix elements $v\leq w$ of $I(d)$. Fix any term order on the
monomials in the standard graded
polynomial ring $P:=k[X_\beta\st\beta\in\roots^v]$ such
that, for any good admissible pair $\ap$,  the initial term of
$f_\ap$ 
is $\mon_\theta$. 
Then the\/
$f_\ap$ as $\ap$ varies over all good admissible pairs
form a Gr\"obner basis
with respect to this term order for the ideal
$(f_\ap\st
\ap\text{ admissible pair, } \aptop{\ap}\nleq w)$
defining the
tangent cone to the Schubert variety $X_w$ at the point $e^v$.\eprop
\bproof   
Denoting by $I$ the ideal
$(f_\ap\st\ap\text{ admissible pair,
$\aptop{\ap}\nleq w$})$ defining the tangent cone, 
by $\init{f}$ the
initial term of a polynomial $f$ in the fixed term order,  and by $\init{I}$
the ideal $\left(\init{f}\st f\in I\right)$,  we clearly have a graded
surjection
\[
P/
\left(
\init{f_\ap}\st\ap\text{ good admissible pair,
$\aptop{\ap}\nleq w$}\right) 
\surjection P/\init{I}
\]
The assertion is that this map is an isomorphism. To 
prove this,  it is enough to show that the Hilbert function
of the quotient ring dominates that of the ring on the left.

The Hilbert function of $P/\init{I}$ is the same as that of the tangent cone
$P/I$,  and so by Theorem~\ref{tmain} its value at a positive
integer $m$ is the cardinality of
the set $\svwm$ of $w$-dominated
monomials in $\roots^v$ of degree $m$.  It therefore suffices to
observe that a monomial
in $\roots^v$ that is not $w$-dominated belongs to the ideal
$\left(
\init{f_\ap}\st\ap\text{ good admissible pair,
$\aptop{\ap}\nleq w$}\right)$.
Given such a monomial, choose a $v$-chain
$\beta_1>\ldots>\beta_t$ in it
such that $w\ngeq s_{\beta_1}\cdots s_{\beta_t}v$.  Applying the lemma below
to the this $v$-chain,  we get a good admissible pair $\ap$ such that
the initial term of $f_\ap$ in the fixed term order is
$\{\beta_1,\ldots,\beta_t\}$.\eproof
\blemma\label{lgrob}
Fix elements $v\leq w$ of $I(d)$.
Let 
$\beta_1>\ldots>\beta_r$ be a $v$-chain of elements in $\pos^v$
such that $w\ngeq s_{\beta_1}\cdots s_{\beta_r}v$.
The element of\/
${\rm SM}^{v,v}$\/ that corresponds under the bijection of
Proposition~\ref{pbijection} to the element 
$\{\beta_1,\ldots,\beta_r\}$ of\/ $T^v$
consists of a single admissible pair $\ap=(t,u)$ which is good and for which
$\mon_\theta=\{\beta_1,\ldots,\beta_r\}$.
\elemma

\bproof
We follow the notation of the statement and
proof of Proposition~\ref{pbijection}.    
Set $\mon =\{\beta_1,\ldots,\beta_r\}
\cup\{\beta_1^\hash,\ldots,\beta_r^\hash\}$,   $\pi(\mon)=(t,\mon')$,
and $\pi(\mon')=(u,\min)$.      We know from that
proposition that $(t,u)$ is an admissible pair and that
$v\leq u$.    Since the bijection ${\rm SM}^{v,v}\simeq T^v$ respects
domination,   and $s_{\beta_1}\cdots s_{\beta_r}v\nleq w$ by hypothesis,
it follows that $t\nleq w$.    It remains to be seen that
$\min$ is the empty monomial and
that $\mon_\theta=\{\beta_1,\ldots,\beta_r\}$.  These will follow from 
explicit computation which we now perform.

Set $\beta_j=(r_j,c_j)$;  then
$\beta_j^\hash=(c_j^*,r_j^*)$. 
Since $\beta_j\in\pos^v$,  we have $c_j\leq r_j^*$ and $c_j<r_j$,  and so
$\{c_1,\ldots,c_t\}\subseteq [d]$. 

Consider the two-step partitioning of $\mon$ into the $\mon_j$
and the $\mon_j$ into blocks.   It is readily seen that
$\mon_j=\{\beta_j,\beta_j^\hash\}$.
Let $s$ be the largest integer, $0\leq s\leq r$, 
such that $r_s\geq r_s^*$---set
$s=0$ if $r_1<r_1^*$.    In other words,
$s$ is the least integer, $0\leq s\leq t$, such that 
$\{r_{s+1},\ldots,r_t\}\subseteq [d]$. 
For $1\leq j\leq s$,   $\mon_j$ is a single block.   For $s<j\leq t$,
the block decomposition of $\mon_j$ is
$\{\beta_j\}\cup\{\beta_j^\hash\}$.

The following expressions for $\mon'$, $u$, $t$, and $\min$
follow readily from 
the definition of the map $\pi$:
$$
\mon'  =  \{(r_1,r_1^*),\ldots,(r_t,r_t^*)\} \quad {\rm and} \quad 
u  =  \left( v\cup\{r_{1},\ldots,r_s\}\right)
	\setminus \left( \{r_{1}^*,\ldots,r_s^*\}\right),
$$
and 
\[ %\begin{array}{c}
t  =  \left( v\cup\{c_1^*,\ldots,c_t^*\}\cup\{r_{s+1},\ldots,r_t\} \right)
  \setminus \left( \{c_1,\ldots,c_t\}\cup\{r_{s+1}^*,\ldots,r_t^*\} \right),
									\]
and $\min$ is the empty monomial.   Now, upon letting $[d]^{\rm c}:=\{d+1,\ldots,2d\}$, put
$$
\theta :=  \left(t\cap [d] \right)\cup \left( u\cap [d]^{\rm c} \right) = 
\left( v\cup\{r_1,\ldots,r_t\}\right)	\setminus\{c_1,\ldots,c_t\}
$$
and
$$
\theta^\hash :=  \left(u \cap [d] \right)\cup \left( t \cap \right [d]^{\rm c}) = 
\left(v\cup\{c_1^*,\ldots,c_t^*\}\right) \setminus\{r_1^*,\ldots,r_t^*\}.
$$
Immediately from the definition of $\mon_\theta$,  we have
\[ 
\mon_\theta=\{(r_1,c_1),\ldots,(r_t,c_t)\}\quad\text{and}\quad
 \mon_{\theta^\hash}=\{(c_1^*,r_1^*),\ldots,(c_t^*,r_t^*)\}. %\Box 
\]
Thus the lemma is proved.
%\eprooftwo
\eproof

We finish by listing some term orders that satisfy the requirement
of Proposition~\ref{pgroebner}.    Fix notation and
terminology as in \S15.2 of \cite{ebud}.  Defined below are
eight total orders $\poone$ through $\po_8$ 
on the elements of $\pos^v$.  Let $>_j$ denote also any
total order on $\roots^v$ that refines the total order $>_j$ on $\pos^v$.
The homogeneous
lexicographic orders induced by $\poone$,  $\potwo$, $\po_7$, $\po_8$
and the
reverse lexicographic orders induced by $\po_4$, $\po_6$ 
satisfy the requirement of Proposition~\ref{pgroebner}.   The reverse
lexicographic orders induced by $\po_3$, $\po_5$ do not satisfy the requirement.

The total orders $\poone$ through $\po_8$ on $\pos^v$ are as follows:
	for $(r,c)$ in $\pos^v$ and $(a,b)$ in $\roots^v\setminus\pos^v$,
	\begin{enumerate}
	\item	$(r,c)\poone(r',c')$ if either (a)~$r<r'$ or
			(b)~$r=r'$ and $c>c'$.
	\item	$(r,c)\potwo(r',c')$ if either (a)~$c>c'$ or
			(b)~$c=c'$ and $r<r'$.
	\item	$(r,c)\pothree(r',c')$ if either (a)~$r<r'$ or
			(b)~$r=r'$ and $c<c'$.
	\item	$(r,c)\pofour(r',c')$ if either (a)~$c>c'$ or
			(b)~$c=c'$ and $r>r'$.
	\item	$(r,c)\po_5(r',c')$ if either (a)~$c<c'$ or
			(b)~$c=c'$ and $r<r'$.
	\item	$(r,c)\po_6(r',c')$ if either (a)~$r>r'$ or
			(b)~$r=r'$ and $c>c'$.
	\item	$(r,c)\po_7(r',c')$ if either (a)~$c<c'$ or
			(b)~$c=c'$ and $r>r'$.
	\item	$(r,c)\po_8(r',c')$ if either (a)~$r>r'$ or
			(b)~$r=r'$ and $c<c'$.
			\end{enumerate}

\section*{Acknowledgments}
\label{sacknow}
A part of this work was done when the first named author visited the Institute of
Mathematical Sciences, Chennai, during March 2003, and Purdue University, West Lafayette,
during Spring 2004. He would like to thank both these institutions for their hospitality
and support. Also, he gratefully acknowledges the travel support from the 
Commission on Development and Exchanges of the International Mathematical Union for his 
visit to Purdue. 

\bibliographystyle{amsplain}

\end{document}